\documentclass[11pt]{article}

\usepackage[utf8]{inputenc}

\usepackage[square,numbers]{natbib}
\usepackage{fullpage}
\usepackage{xcolor}
\usepackage{amsmath}
\usepackage{amssymb}
\usepackage{bold-extra}
\usepackage{amsthm}
\usepackage{mathtools}
\usepackage{setspace}
\usepackage{enumitem}
\newcommand{\Ab}{\mathbb{A}}
\newcommand{\Eb}{\mathbb{E}}
\newcommand{\Hb}{\mathbb{H}}
\newcommand{\Kb}{\mathbb{K}}

\newcommand{\Pb}{\mathbb{P}}
\newcommand{\Qb}{\mathbb{Q}}
\newcommand{\Rb}{\mathbb{R}}
\newcommand{\Sb}{\mathbb{S}}

\newcommand{\Xb}{\mathbb{X}}
\newcommand{\Yb}{\mathbb{Y}}

\newcommand{\Bc}{\mathcal{B}}

\newcommand{\Xc}{\mathcal{X}}
\newcommand{\Hc}{\mathcal{H}}
\newcommand{\Tc}{\mathcal{T}}

\renewcommand{\bar}{\overline}
\renewcommand{\vec}{\mathbf}
\newcommand{\hto}{\xrightarrow{d_H}}

\newcommand{\inv}{^{-1}}

\newcommand{\oneAlt}[1]{\mathbf{1}{#1}}

\newcommand{\defeq}{:=}

\newcommand{\pomdpOne}{{POMDP}${}_{1}$}
\newcommand{\pomdpTwo}{{POMDP}${}_{2}$}
\newcommand{\comdp}{{COMDP}}

\DeclareMathOperator{\Uniform}{\mathrm{Uniform}}

\DeclareMathOperator{\diag}{\mathbf{Diag}}

\DeclarePairedDelimiter{\set}{\{}{\}}
\DeclarePairedDelimiter{\abs}{|}{|}
\DeclarePairedDelimiter{\norm}{\|}{\|}

\newtheorem{lemma}{Lemma}[section]
\newtheorem{theorem}[lemma]{Theorem}
\newtheorem{corollary}[lemma]{Corollary}
\newtheorem{definition}[lemma]{Definition}

\newtheorem{example}[lemma]{Example}

\newtheorem{remark}[lemma]{Remark}

\newtheorem{assumption}[lemma]{Assumption}

\newtheorem{mainCondition}{Diffeomorphic Condition}
\newtheorem{unAssumption}{Assumption}

\allowdisplaybreaks

\title{%
    Continuity of Filters for Discrete-Time Control Problems Defined by Explicit Equations\footnote{
       This research was partially supported by  SUNY System Administration under  SUNY Research Seed Grant Award 231087 and by the U.S. Office of Naval Research (ONR) under Grants N000142412608 and N000142412646.
    }}

\author{
    Eugene A. Feinberg\thanks{%
        Department of Applied Mathematics and Statistics, State University of New York at Stony Brook, Stony Brook, NY 11794-3600, USA, (email: {eugene.feinberg\@stonybrook.edu}).
    } \and
    Sayaka Ishizawa\thanks{%
        Department of Applied Mathematics and Statistics, State University of New York at Stony Brook, Stony Brook, NY 11794-3600, USA, (e-mail: {sayaka.ishizawa\@alumni.stonybrook.edu}).
    } \and %
    Pavlo O. Kasyanov\thanks{%
        Institute for Applied System Analysis, National Technical University of Ukraine “Igor Sikorsky Kyiv Polytechnic Institute", Peremogy ave., 37, build, 35, 03056, Kyiv, Ukraine, (e-mail: {kasyanov\@i.ua}).
    } \and
    David N. Kraemer\thanks{%
        Department of Applied Mathematics and Statistics, State University of New York at Stony Brook, Stony Brook, NY 11794-3600, USA, (e-mail: {david.n.kraemer\@gmail.com}).
    }
}

\date{\today}

\begin{document}


\maketitle

\begin{abstract}
    Discrete time control systems whose dynamics and observations are described by stochastic equations are common in engineering, operations research, health care, and economics. For example, stochastic filtering problems are usually defined via stochastic equations. These problems can be reduced to Markov decision processes (MDPs) whose states are posterior state distributions, and transition probabilities for such MDPs are sometimes called filters.  This paper investigates sufficient conditions on transition and observation functions for the original problems to guarantee weak continuity of  the filter.  Under mild conditions on cost functions, weak continuity implies the existence of optimal policies minimizing the expected total costs, the validity of optimality equations, and convergence of value iterations to optimal values.  This paper uses recent results on weak continuity of filters for partially observable MDPs defined by transition and observation probabilities. It develops a criterion of weak continuity of transition probabilities and a sufficient condition for continuity in total variation of transition probabilities. The results are illustrated with applications to filtering problems.
\end{abstract}



\section{Introduction} \label{sec:introduction}

This paper provides conditions for weak continuity of filters for control systems with incomplete information defined by stochastic equations. A filter is the transition probability for a completely observable controlled process, often called a belief Markov Decision Process (MDP) or a Completely Observable MDP (COMDP), whose states are posterior probabilities of the states of the original process. 
Weak continuity of a filter and appropriate continuity conditions on one-step costs imply the existence of optimal policies, validity of optimality equations, and convergence of value iterations for problems with expected total costs.   The  analysis in this paper is based on recently developed in \cite{feinberg_partially_2016, feinberg_markov_2022, feinberg_2023_semi, kara_2019_weak} sufficient conditions for weak continuity of filters for Partially Observable MDPs (POMDPs).  POMDPs are stochastic control models with incomplete information defined by transition and observation probabilities. 

Thus there are two approaches for modeling control problems with control information: (a) to use stochastic equations, as this is done in stochastic filtering, or (b) to use transition and observation probabilities also called transition and observation kernels, as this is done for POMDPs.   The stochastic equation approach and the POMDP approach are equivalent; see Aumann's lemma, which is Theorem~\ref{thm:kernels:aumann} below. Aumann's lemma links transition probabilities and the corresponding transition functions.  To investigate weak continuity of filters for problems defined by stochastic equations, this paper studies relations between properties of transition functions and continuity of the corresponding transition probabilities.

The natural way to find optimal policies for POMDPs is to reduce the initial problem to a similar problem for a COMDP.  If an optimal policy for a COMDP is found, then it can be used to construct an optimal policy for the initial problem that minimizes the expected total discounted costs~\cite{bertsekas_stochastic_1996,dynkin_1979_controlled,hernandez-lerma_adaptive_1989}.~\citet{aoki_1965},~\citet{astrom_1965},~\citet{dynkin_1965}, and~\citet{shiryaev_1966} introduced this approach for problems with finite state, observation, and action spaces.~\citet{sawaragi_1970} extended it to problems with countable state spaces, and~\citet{rhenius_1974} and~\citet{yushkevich_1976_reduction} to Borel state spaces. For linear systems with  Gaussian noises, the distributions of states are Gaussian, and they are defined by mean vectors and covariance matrices.  In this case, the Kalman filter completely characterizes  transition probabilities.

In general, the existence of optimal policies for POMDPs is a nontrivial question even for problems with expected total costs.  For MDPs in general and for COMDPs in particular, for expected total costs the existence of optimal policies and the validity of other important properties, including the existence of solutions to optimality equations and convergence of value iterations, follow from continuity properties of transition probabilities and of  one-step costs.  For MDPs with compact action sets, the corresponding properties are weak continuity of transition kernels and lower semicontinuity of one-step costs; \citet{schal_1972_continuous,schal_1993_average}. Without the compactness assumption, the corresponding properties are weak continuity of transition kernels and $\Kb$-inf-compactness of one-step costs; \citet{feinberg_average_cost_2012,feinberg_partially_2016}. There is also a parallel theory for MDPs with setwise continuous transition probabilities  \cite{feinberg_2021_mdps,schal_1972_continuous,schal_1993_average}, but so far it has not found applications to problems with incomplete observations.

While lower semicontinuity and $\Kb$-inf-compactness of costs are preserved by the reduction of a POMDP to a COMDP, weak continuity of the transition probability may not be preserved \cite[Theorem 3.3, Lemma 2.1, and Examples 4.2 and 4.3]{feinberg_partially_2016}. For a long time, sufficient conditions for weak continuity of filters were unavailable.~\citet[pp. 90-93]{hernandez-lerma_adaptive_1989} and~\citet[Section 2]{runggaldier_1994_approximations} introduced some particular conditions assuming, among other conditions, weak continuity of the transition kernel  and continuity in total variation of the observation kernel.~\citet{feinberg_partially_2016} proved that these two conditions imply weak continuity of filters.~\citet{kara_2019_weak} gave another proof of this fact and proved that continuity of the transition kernel is sufficient if the observation kernel  does not depend on observations.~\citet{feinberg_markov_2022,feinberg_2023_semi} introduced the notion of semi-uniform Feller transition probabilities and showed that this property is preserved when a problem with incomplete observations is reduced to a COMDP, and this is also true for more general problems with incomplete information than POMDPs. This property implies weak continuity of the filter.  In particular, these facts provide another proof that weak continuity of the transition kernel and continuity of the observation kernel in total variation for a POMDP imply weak continuity of the filter.  They also imply that, if the transition kernel  is continuous in total variation and the observation kernel  is continuous in total variation in the control parameter, then the filter is weakly continuous.

This paper provides sufficient conditions 
on the transition and observation stochastic equations defining the dynamics of the system for weak continuity of the filter.  Such conditions were considered by~\citet[Section 8.1]{feinberg_partially_2016} for problems with real-valued states, observations, and noises. This paper obtains stronger results for  problems with Euclidean state, observation, and noise spaces.  The important technical result is Theorem~\ref{thm:kernels:total_variation}(b) providing a sufficient condition for continuity  of a transition probability in total variation. The proof of this theorem uses continuity of Diffeomorphic images of a compact set in $\Rb^n$ both in Hausdorff and  Fr\'echet-Nikodym senses stated in Lemma~\ref{lem:proofs:images_converge_in_measure}. A relevant result is provided in \cite{KY}.  By using Scheff\'e's theorem, \cite[Example~1.3.2(v)]{KY}  provides a simple proof of continuity of transition probability in total variation under two additional assumptions: (i) the noise density is continuous, and (ii) parametric images of the sample space are constant, and assumption (ii) is not stated in \cite{KY}; see Remark~\ref{rem:Shaffe} below for details.

This paper is organized as follows. Section~\ref{sec:notation} introduces the notations and definitions used throughout this paper. Section~\ref{sec:filters} formulates the main problem, introduces assumptions, and states Theorem~\ref{thm:filters:main_result}, which is the main result of this paper. Section~\ref{sec:feller} provides sufficient conditions for semi-uniform Feller continuity of the joint distribution of the state-observation pair at the next state given the current state-action pair. Section~\ref{sec:kernels} explains the equivalence between processes defined by stochastic equations and by stochastic kernels (Theorem~\ref{thm:kernels:aumann} introduced by~\citet[Lemma F]{aumann_1964_mixed} and its Corollary~\ref{cor:kernels:aumann_multidim}) and provides criteria and sufficient conditions for weak continuity and for continuity in total variation of stochastic kernels defined by equations (Theorems~\ref{thm:kernels:weak} and~\ref{thm:kernels:total_variation}). Section~\ref{sec:pomdps} formulates the optimization problem, optimality equations, and provides sufficient conditions for the existence of optimal policies   and convergence of value iterations for problems defined by stochastic equations.  Section~\ref{sec:applications} considers applications to stochastic filtering, and it focuses on state space models with additive and multiplicative transition and observation noises, including linear state space models and inventory control models. Section~\ref{sec:pomdp1s} describes a parallel theory for a different version of POMDP models,  which has been considered in the operations research literature. Section~\ref{sec:proofs} provides the proofs to the main results of this paper.

\section{Notation and Background Definitions} \label{sec:notation}

For  a metric space $\Sb,$ we denote by $\rho_\Sb$ its metric and denote by $\Bc(\Sb)$ its Borel $\sigma$-algebra; that is,  $\Bc(\Sb)$  is the $\sigma$-algebra generated by the open subsets of $\Sb.$ For $S \in \Bc(\Sb)$ we usually consider the metric space $(S,\rho_S),$ where $\rho_S(s_1,s_2)= \rho_\Sb(s_1,s_2)$ for $s_1,s_2\in S.$ For the Borel $\sigma$-algebra $\Bc(S)$ on $S,$ the equality $\Bc(S)=\set{B\in \Bc(\Sb): B\subset S}$ holds.
For   $f:S\to E,$ where    $S\subset\Sb,$    $E\subset\mathbb{E},$ and $\Sb$ and $\mathbb{E}$ are metric spaces, we always view $S$ and $E$ as metric spaces.  In particular, 
for $f:\Sb\to\Eb,$  continuity of $f:S\to E$ means than the restriction of $f$ on $S$ is a continuous mapping with values in $E.$ For metric spaces $\Sb_1$ and $\Sb_2$ we consider the metric space $\Sb_1\times\Sb_2$ with the product metric
$\rho_{{\Sb_1\times\Sb_2}}((s_1',s_2'),(s_1'',s_2'')):=[\rho^2_{\Sb_1}(s_1',s_1'')+\rho^2_{\Sb_2}(s_2',s_2'')]^{\frac{1}{2}},$ $(s_1',s_2'),(s_1'',s_2'')\in \Sb_1\times\Sb_2.$

The set of probability measures on the measurable space $(\Sb, \Bc(\Sb))$ is denoted by $\Pb(\Sb).$ The set $\Pb(\Sb)$ is endowed with the topology of weak convergence of probability measures, that is, $p^{(k)} \to p$ in $\Pb(\Sb)$ if $\int_{\Sb} f(s) \ p^{(k)}(ds) \to \int_{\Sb} f(s) \ p(ds)$ for all bounded continuous functions $f : \Sb \to \Rb.$ Let $\Sb_1$  and $\Sb_2$ be Borel spaces; that is, Borel subsets of Polish (separable, complete metric) spaces. We recall that a stochastic kernel on $\Sb_1$ given $\Sb_2$ is a function $\kappa : \Bc(\Sb_1) \times \Sb_2 \to [0,1],$ written $\kappa(B|s_2),$ such that
\begin{enumerate}[label=(\alph*)]
    \item for each $B \in \Bc(\Sb_1),$ the map $s_2 \mapsto \kappa(B|s_2)$ is a Borel measurable function,
    \item for each $s_2 \in \Sb_2,$ the map $B \mapsto \kappa(B|s_2)$ is a Borel probability measure on $\Sb_1.$
\end{enumerate}
A stochastic kernel $\kappa$ on $\Sb_1$ given $\Sb_2$ is weakly continuous if the probability measures $\kappa(\:\cdot\:|s_2')$ converge weakly to $\kappa(\:\cdot\:|s_2)$ in $\Pb(\Sb_1)$ as $s_2' \to s_2,$ and continuous in total variation if $\sup_{B \in \Bc(\Sb_1)} \abs{\kappa(B|s_2') - \kappa(B|s_2)} \to 0$ as $s_2' \to s_2.$ Continuity in total variation implies weak continuity.

In this paper, the variables $d,$ $m,$ and $n$ always refer to positive integers. Given two measures $p_1, p_2$ {  on the same measurable space} we write $p_1 \ll p_2$ if $p_1$ is absolutely continuous with respect to $p_2.$ 
Lebesgue measure on $\Rb^n$ is denoted $\lambda^{[n]}.$ Let $D_x g = \frac{\partial g}{\partial x}$ denote the Jacobian of  a differentiable function $g : \Rb^n \to \Rb^n.$ When we consider the sufficiently smooth function  $g(s_2,x)$ on $\mathbb{S}_2\times\Rb^n,$ 
we denote its Jacobian in $x$ either by $D_x \phi(s_2,x)$ or by $D_x \phi_{s_2}(x).$

A function $f : \Sb \to [-\infty, +\infty]$ is lower semicontinuous at $s \in \Sb$ if $\liminf_{s' \to s} f(s') \geq f(s).$ If $f$ is lower semicontinuous at each $s \in \Sb,$ then $f$ is lower semicontinuous on $\Sb.$ The function $f$ is inf-compact if for all $\gamma \in \Rb,$ the level set $\set{s \in \Sb: f(s) \leq \gamma}$ is compact. A function $f : \Sb_1 \times \Sb_2 \to [-\infty,+\infty]$ is $\Kb$-inf-compact on $\Sb_1 \times \Sb_2$ if for all nonempty compact sets $C \subset \Sb_1$ and  for all $\gamma \in \Rb,$ the sublevel sets $\set{(s_1,s_2) \in C \times \Sb_2\, :\, f(s_1,s_2) \leq \gamma}$ are compact.

We denote by $\Xb,$ $\Yb,$ and $\Ab$ the state, observation, and action spaces, respectively. We denote by $\Xc$ and $\Hc$  the spaces of  state and observation noises.  In general, $\Xb,$ $\Yb,$   $\Ab,$ $\Xc,$ and $\Hc$ are assumed to be Borel spaces, that is, they are Borel subsets of complete separable metric spaces. 
In some results there are additional assumptions that some of these spaces are Euclidian.  The discrete time parameter is $t=0,1,\ldots,$ and $x_t,$ $y_t,$ $a_t,$ $\xi_t,$ and $\eta_t$ denote the state, observation, control, state noise, and observation noise at time $t$ respectively, where  $x_t\in \Xb,$ $y_t\in \Yb,$   $a_t\in \Ab,$ $\xi_t\in \Xc,$ and $\eta_t\in \Hc.$

For a Borel space $\Omega,$ for   a Borel measurable function $\phi : \Sb_2 \times \Omega \to \Sb_1,$ and for a probability measure $p \in \Pb(\Omega),$ let us define the stochastic kernel $\kappa$ on $\Sb_1$ given $\Sb_2,$
\begin{align} \label{eq:notation:stochastic_kernel}
    \kappa(B|s_2) &\defeq \int_{\Omega} \mathbf{1}\{\phi(s_2,\omega) \in B\} \ p(d\omega), &&
    B \in \Bc(\Sb_1), \quad
    s_2 \in \Sb_2.
\end{align}
The following properties of $\phi$ will be used throughout the paper.

\begin{definition}[Continuity in distribution, total variation, and probability]\label{def:notation:phi_continuity}
    Let  $\Sb_1,$ $\Sb_2,$ and $\Omega$ be Borel spaces, and let  $p \in \Pb(\Omega)$ be a probability measure on $(\Omega, \Bc(\Omega)).$ A Borel function $\phi : \Sb_2 \times \Omega \to \Sb_1$ is continuous
    \begin{enumerate}
        \item[{\rm(i)}] \label{item:notation:continuity_in_distribution} in distribution $p$ (weakly continuous) if the function $s_2 \mapsto \int_\Omega f(\phi(s_2,\omega))\ p(d\omega)$ is continuous on $\Sb_2$ for every bounded continuous function $f : \Sb_1 \to \Rb;$
        \item[{\rm(ii)}] \label{item:notation:continuity_in_total_variation} in total variation with respect to ({\it wrt}) $p$ if for each $s_2 \in \Sb_2,$
        \begin{equation}
            \lim_{s_2' \to s_2}
            \sup_{B \in \Bc(\Sb_1)} \left|
                \int_{\Omega} \mathbf{1}\{\phi(s_2', \omega) \in B\} \
                -
                \mathbf{1}\{\phi(s_2, \omega) \in B\} \ p(d\omega)
            \right| = 0;
           \end{equation}
        \item[{\rm(iii)}] \label{item:notation:continuity_in_probability} in probability $p$ if $\phi(s_2',\:\cdot\:) \xrightarrow{p} \phi(s_2,\:\cdot\:)$ as $s_2' \to s_2$ for each $s_2\in \Sb_2,$ that is, for each $s_2 \in \Sb_2$ and each $\varepsilon > 0,$
        \begin{equation}
            \lim_{s'_2 \to s_2} p(\set{
                \omega\in \Omega \,:\, \rho_{\Sb_1}(\phi(s_2',\omega),\phi(s_2,\omega)) \ge \varepsilon
            }) = 0.
        \end{equation}

    \end{enumerate}
\end{definition}
 We note that continuity in probability is called stochastic continuity~\cite[p. 30]{wentzell_1996_course}.

If $\phi$ is continuous in total variation {\it wrt} $p,$ then it is stochastically continuous. It is well-known that continuity in probability implies continuity in distribution, while the opposite statement is false~\cite[p. 30]{wentzell_1996_course}.

For the Euclidean space $\Rb^n,$ its standard norm is denoted by $\norm{\cdot}.$ For a function {$\phi : \Sb_2 \times \Rb^n \to \Sb_1,$ we will frequently assume that $\Sb_1 = \Rb^n,$} 
and that the function $\phi$ satisfies the following condition.
\begin{mainCondition} \label{cond:notation:differentiability_condition}
    For the metric space $\Sb_2,$ open set $\Omega \subset \Rb^n,$  and a function $\phi : \Sb_2 \times \Rb^n  \to \Rb^n,$ the following statements hold:
    \begin{enumerate}
        \item[{\rm(i)}] $\phi$ is continuous on $\Sb_2 \times \Omega;$
        \item[{\rm(ii)}] $D_\omega\phi(s_2,\omega)$ exists for all $s_2 \in \Sb_2$ and $\omega \in \Omega;$
        \item[{\rm(iii)}] the matrix $D_\omega \phi(s_2,\omega)$ is nonsingular for all $s_2 \in \Sb_2$ and $\omega \in \Omega;$
        \item[{\rm(iv)}] the function $(s_2, \omega) \mapsto D_\omega \phi(s_2,\omega)$ is continuous on $\Sb_2 \times \Omega;$
        \item[{\rm(v)}] for each $s_2 \in \Sb_2$ the function $\omega \mapsto \phi(s_2,\omega)$ is a one-to-one mapping of $\Omega$  onto $\phi(s_2,\Omega).$
    \end{enumerate}
\end{mainCondition}


\begin{remark}\label{rem:new}{\rm  Let us use the notation $\phi_{s_2}(\omega):=\phi(s_2,\omega).$    Diffeomorphic Condition (v) implies that for each $s_2\in\Sb_2$ the mapping $\phi_{s_2}:\Omega\to \phi_{s_2}(\Omega)$ has the inverse mapping  $\phi^{-1}_{s_2}:\phi_{s_2}(\Omega)\to \Omega.$ The Diffeomorphic Condition and the global inverse function theorem imply that for each $s_2\in\Sb_2$ the function  $\phi^{-1}_{s_2}$ is continuously differentiable on $\phi_{s_2}(\Omega).$ Diffeomorphic Condition (iii) and the inverse function rule imply  $D_\omega \phi_{s_2}^{-1}(\omega)\ne 0$ for all all $s_2\in\Sb_2$ and for all $\omega\in \phi_{s_2}(\Omega).$  Let us consider the function $\tilde{\phi}:\Sb_2\times\Omega\to \Sb_2\times\Rb^n,$ where $\tilde{\phi}(s_2,\omega)=(s_2,\phi_{s_2}(\omega)).$ We observe that for each $\tilde{S}_2\subset\Sb_2$ and $\tilde{\Omega}\subset\Rb^n,$ the set $\tilde{\phi}(\tilde{S}_2\times\tilde{\Omega})$ is the graph of the function $\phi:\tilde{S}_2\times\tilde{\Omega}\to\Rb^n.$  Let $\Kb\subset \Omega$ be nonempty and compact. Then the mappings $\phi^{-1}_{s_2}(\omega):\tilde{\phi}(\Sb_2\times\Kb)\to\Kb$
and $\det D_\omega \phi^{-1}_{s_2}(\omega):\tilde{\phi}(\Sb_2\times\Kb)\to\Rb$ are continuous. This is true because of the following reasons.  Let us consider a sequence $\{s_2^{(k)}\}_{k=1}^\infty$ in $\Sb_2$ convergent to a point $\tilde{s}_2\in \Sb_2.$  The set $S_2:=\{\tilde{s}_2,  s_2^{(1)},  s_2^{(2)},\cdots\}$ is a compact subset of $\Sb_2.$ Thus, $S_2\times\Kb$ is a compact subset of $\Sb_2\times\Kb.$  
The function $\tilde{\phi}:\Sb_2\times\Omega\to \Sb_2\times\Rb^n$  is continuous as a composition of two continuous functions. Diffeomorphic Condition (v) implies that this function is bijective.
 In view of \cite[Theorem~26.6]{munkrestopology},  the mapping $\tilde{\phi}:S_2\times\Kb\to \tilde{\phi}(S_2\times \Kb)$ is a homeomorphism.  Therefore, the mapping $\tilde{\phi}^{-1}:\tilde{\phi}(S_2\times \Kb)\to \Sb_2\times\Kb$ is  continuous, where $\tilde{\phi}^{-1}(s_2,\omega)=(s_2, \phi^{-1}_{s_2}(\omega))$ for $(s_2,\omega)\in \tilde{\phi}(S_2\times \Kb).$ 
 Since $\{s_2^{(k)}\}_{k=1}^\infty$ is an arbitrary sequence in $\Sb_2$ convergent to an arbitrary $\tilde{s}\in\Sb_2,$ the mapping $\tilde{\phi}^{-1}:\tilde{\phi}(\Sb_2\times\Kb)\to\Kb$ is continuous.  Therefore, the mapping $\phi^{-1}_{s_2}(\omega):\tilde{\phi}(\Sb_2\times\Kb)\to\Kb$
  is continuous. In addition, the inverse function rule and the assumption $D_\omega \phi_{s_2}(\omega)\ne 0$ imply that the function $\det D_\omega \phi^{-1}_{s_2}(\omega):\tilde{\phi}(\Sb_2\times\Kb)\to\Rb$ is continuous.}
\end{remark}

We also note that, in the one-dimensional { case} 
 Diffeomorphic Condition~(v) follows from Diffeomorphic Conditions~(i-iv) since each function $\phi_{s_2}(\omega)$ is continuous and strictly monotonic. If $n>1,$ finding broad sufficient conditions for injectivity of the mapping $\omega\mapsto \phi_{s_2}(\omega)$ over the entire set $\Omega$ is a challenging mathematical problem. 
However, as shown in Section~\ref{sec:applications} below, in the case of 
additive or multiplicative noises, Diffeomorphic Condition~(v) is satisfied  if Diffeomorphic Conditions~(i-iv) are satisfied.

\section{Continuity of Filter Kernels Generated by Stochastic Equations} \label{sec:filters}

Let
$\Xb,$ $\Yb,$ and $\Ab$ be Borel spaces of states, observations, and actions, respectively, and $\Xc, \Hc$ be Borel spaces of state and observation noises.
We consider a discrete-time control system with dynamics and observations defined for time $t = 0,1, \dots$ by stochastic equations
\begin{subequations}\label{eq:filters:model}
    \begin{align}
        x_{t+1} &= F(x_t, a_t, \xi_t), && x_t \in \Xb, \quad a_t \in \Ab, \quad \xi_t \in \Xc, \label{eq:filters:transitions} \\
        y_{t+1} &= G(a_t, x_{t+1}, \eta_{t+1}), && a_t \in \Ab, \quad x_{t+1} \in \Xb, 
        \quad \eta_{t+1} \in \Hc, \label{eq:filters:observations}
    \end{align}
\end{subequations}
where $F : \Xb \times \Ab \times \Xc \to \Xb$ is the Borel-measurable functions representing the transition dynamics of the system, and $G : \Ab \times \Xb \times \Hc \to \Yb$ is the Borel-measurable function representing the observations. In addition, there is the probability distribution $p_0 \in \Pb(\Xb)$ of the initial state $x_0,$ and $y_0= G_0(x_0,\eta_0),$  where $G_0 : \Xb \times \Hc \to \Yb$ a Borel-measurable function defining the initial observation, and
$\set{\xi_t}_{t=0}^{\infty}$ and $\set{\eta_t}_{t=0}^{\infty}$ are sequences of independent and identically distributed (iid) random variables with distributions $\mu \in \Pb(\Xc)$ and $\nu \in \Pb(\Hc),$ respectively. It is also assumed that the sequence $\set{x_0, \xi_0, \eta_0, \xi_1, \eta_1, \dots}$ is mutually independent. The process evolves as follows.  The initial hidden state $x_0$ has the distribution $p_0,$ and the initial observation is $y_0 = G(x_0, \eta_0).$  At each epoch $t = 0, 1, \dots,$  the decision maker, who  knows the initial state distribution $p_0,$ the previous and current observations, and the previously selected actions $p_0,y_0,a_0,\ldots,y_{t-1},a_{t-1},y_t,$ selects the action $a_t \in \Ab.$ The next state $x_{t+1}$ and observation $y_{t+1}$ are determined by equations \eqref{eq:filters:model}. 

The evolution of a system defined by stochastic equations such as~\eqref{eq:filters:model} can also be represented in terms of stochastic kernels, and the corresponding model is a POMDP.  For POMDPs, transitions of states are defined by a stochastic kernel $\Tc$ on $\Xb$ given $\Xb\times\Ab$ called the transition kernel or transition probabilities, and observations are defined by a stochastic kernel $Q$ on $\Yb$ given $\Ab\times\Xb$ called the observation kernel or observation probabilities.  For POMDPs, the state $x_{t+1}$ is defined by the distribution $\Tc(\:\cdot\:|x_t,a_t)$ instead of by equality~\eqref{eq:filters:transitions}, and the observation $y_{t+1}$ is defined by the distribution $Q(\:\cdot\:|a_t,x_{t+1})$  instead of by equality~\eqref{eq:filters:observations}. In addition, the initial distribution of observation $y_0$ is defined by a stochastic kernel $Q_0$ on $Y$ given $X,$ that is, $y_0$ is defined by the distribution $Q(\:\cdot\:|x_0);$ see \cite[Chapter 4]{hernandez-lerma_adaptive_1989} or \cite{feinberg_partially_2016} for details.

The model definitions based on stochastic equations and on POMDPs are equivalent.  Indeed, for the functions $F$ and $G$ from~\eqref{eq:filters:model},
\begin{subequations} \label{eq:filters:model_kernels}
    \begin{align}
        \Tc(B|x, a) &= \int_{\Xc} \mathbf{1}\{F(x, a, \xi) \in B\} \ \mu(d\xi), && B \in \Bc (\Xb), \quad x \in \Xb, \quad a \in \Ab, \label{eq:filters:transition_kernel} \\
        Q(C|a, x) & =\int_{\Hc}\mathbf{1}\{G(a, x, \eta) \in C\} \ \nu(d\eta), && C \in \Bc (\Yb), \quad a \in \Ab, \quad x \in \Xb, \label{eq:filters:observation_kernel}
    \end{align}
\end{subequations}
are stochastic kernels; see~\citet[Proposition C.2]{hernandez-lerma_adaptive_1989} or 
\citet[p. 190]{bertsekas_stochastic_1996}. Of course, the function $G_0$ also defines the stochastic kernel $ Q_0(C|a, x)  =\int_{\Hc}\mathbf{1}\{G_0( x, \eta) \in C\}\nu(d\eta),$ where $ C \in \Bc (\Yb),$ $a \in \Ab,$ and $x \in \Xb.$
The reduction of a POMDP to the stochastic system defined by equations~\eqref{eq:filters:model}, with $\Xc=[0,1],$ $\Hc=[0,1],$ and $\mu$ and $\eta$ are uniform distributions,  follows from \citet[Lemma F]{aumann_1964_mixed}; see Theorem~\ref{thm:kernels:aumann} below. However, we do not assume here that $\Xc=[0,1],$ $\Hc=[0,1],$ and $\mu$ and $\eta$ are  uniform distributions because sometimes it can be more convenient to consider other distributions; see, e.g.,  Corollary~\ref{cor:kernels:aumann_multidim}.

The classic approach to solving a POMDP is based on its reduction to a completely observable MDP, whose states are probability distributions of the states of the system, sometimes called belief distributions, and this MDP is sometimes called a filter.  Weak continuity of the transition probabilities of the filter, also called weak continuity of the filter, is an important property.  If one-step costs are $\Kb$-inf-compact and bounded below, then weak continuity of the filter implies the existence of optimal policies for problems with expected total costs, the validity of optimality equations, and convergence of value iterations; see Section~\ref{sec:pomdps} for details.

Theorem~\ref{thm:feller:semi-uniform_feller} describes sufficient conditions on the kernels  $\Tc$ and $Q$ for weak continuity of the filter.  The main goal of this paper is  to describe sufficient conditions on the functions $F$ and $G$ for weak continuity of the filter.  Such results were provided in \cite[Section 8.1]{feinberg_partially_2016} for problems with real-valued states, observations, and noises. In this section described stronger results for more general  state, observation, action, and noise spaces.

Let us consider the standard reduction of a POMDP with a transition kernel $\Tc$ and observation kernel $Q$ to a filtered MDP with complete observations and with the transition kernel $q$ to be defined in~\eqref{eq:filters:filter_kernel}. For more details on this reduction, see~\citet[Section~10.3]{bertsekas_stochastic_1996},~\citet[Chapter~8]{dynkin_1979_controlled},~\citet[Chapter~4]{hernandez-lerma_adaptive_1989},~\citet{rhenius_1974}, or~\citet{yushkevich_1976_reduction}. As usual, we introduce the stochastic kernel on $\Xb \times \Yb$ given $\Xb \times \Ab$ defined by the equation
\begin{align}\label{eq:filters:joint_kernel}
    P(B \times C|x,a) &\defeq
    \int_{B} Q(C|a,x') \ \Tc(dx'|x,a), &&
    B \in \Bc(\Xb), \quad
    C \in \Bc(\Yb), \quad
    x \in \Xb, \quad
    a \in \Ab,
\end{align}
which represents the joint probability of the next state-observation pair 
$(x', y')$ given the state-action pair $(x, a).$  
Consider a measure $z \in \Pb(\Xb)$ representing a belief distribution for the current state $x.$ The joint distribution of the next state-observation pair $(x', y')$ conditioned by the measure $z$ and action $a \in \Ab$ is given as the stochastic kernel
\begin{align} \label{eq:filters:joint_belief_kernel}
    R(B \times C | z,a) &\defeq
    \int_{\Xb} P(B \times C | x,a) \ z(dx), &&
    B \in \Bc(\Xb), \quad
    C \in \Bc(\Yb), \quad
    z \in \Pb(\Xb), \quad
    a \in \Ab,
\end{align}
with marginal distribution for $y'$ given by $R'(C|z,a) \defeq R(\Xb \times C | z,a).$ According to~\citet[Proposition 7.27]{bertsekas_stochastic_1996} there exists a stochastic kernel $H$ on $\Xb$ given $\Pb(\Xb) \times \Ab\times\Yb$ 
such that
\begin{align} \label{eq:filters:joint_belief_conditional_representation}
    R(B \times C | z, a) &=
    \int_C H(B|z,a,y) \ R'(dy|z,a), &&
    B \in \Bc(\Xb), \quad
    C \in \Bc(\Yb), \quad
    z \in \Pb(\Xb), \quad
    a \in \Ab.
\end{align}
Alternatively, the stochastic kernel $H$ defines a Borel measurable mapping $(z,a,y') \mapsto H(z,a,y') \in \Pb(\Xb)$ that describes the evolution of the distributions $z$ by $z'= H(z, a, y').$ We remark that the mapping $H$ is a generalization of filters commonly found in applications. In particular, for linear state space models with a Gaussian initial state distribution $p_0$ and with Gaussian noises, the belief $z$ is a Gaussian distribution parameterized by a mean and a covariance. 
The Kalman filter, which computes the next distribution $z'$ from $z,$ $a$ and the next observation $y',$ computes $H$ explicitly in this case.  Transition dynamics from the current state distribution $z$ to the next state distribution $z'$ is defined by the stochastic kernel $q$ defined as
\begin{align} \label{eq:filters:filter_kernel}
    q(D |z, a) &\defeq
   \int_{\Yb} \mathbf{1}\{H(z, a, y) \in D\} \ R^{\prime}(dy |z, a), &&
    D \in \Bc(\Pb(\Xb)), \quad
    z \in \Pb(\Xb), \quad
    a \in \Ab;
\end{align}
see, e.g., \cite{bertsekas_stochastic_1996, dynkin_1979_controlled, feinberg_partially_2016, hernandez-lerma_adaptive_1989, rhenius_1974, yushkevich_1976_reduction}.
The stochastic kernel $q$ is the transition probability for the filter.

Let us consider again the model defined by equations \eqref{eq:filters:model}.  In view of~\eqref{eq:filters:model_kernels},
\begin{equation} \label{eq:filters:joint_kernel_representation}
    P(B \times C | x,a) =\int_{\Xc}\int_{\Hc} \mathbf{1}\{G(a,F(x,a,\xi),\eta) \in C\} \mathbf{1}\{F(x,a,\xi)\in B\} \ \nu(d\eta) \ \mu(d\xi).
\end{equation}
We observe that formula~\eqref{eq:filters:joint_kernel_representation} employs the assumption of joint independence of the noise variables $\xi$ and $\nu$ mentioned above.

Let us introduce Assumptions~\ref{ass:filters:assumption_1} and \ref{ass:filters:assumption_2}, each of which implies weak continuity of the stochastic kernel~$q.$
\begin{assumption} \label{ass:filters:assumption_1}
    The following statements hold:
    \begin{enumerate}[label=\alph*)]
        \item[\rm(i)] \label{item:filters:F1} the function $((x,a), \xi) \mapsto F(x,a,\xi)$ is continuous in distribution $\mu;$
        \item[\rm(ii)] \label{item:filters:G1} $\Yb=\Rb^n,$ $\Hc$ is an open subset of $\Rb^m,$ $\nu \ll \lambda^{[m]},$ and the function $((a,x),\eta) \mapsto G(a,x,\eta)$ satisfies the Diffeomorphic Condition.
    \end{enumerate}
\end{assumption}

As shown in Section~\ref{sec:kernels}, below, Assumption~\ref{ass:filters:assumption_1}(i) implies weak continuity of the transition kernel $\Tc,$ and Assumption~\ref{ass:filters:assumption_1}(ii) implies continuity in total variation of the observation kernel $Q.$
Assumption~\ref{ass:filters:assumption_1} is the generalization of the assumptions in~\citet[Section 8.1]{feinberg_partially_2016} for $\Xb = \Yb = \Ab = \Xc = \Hc = \Rb$ and $\nu = \Uniform(0,1)$ to finite dimensional spaces. We remark that the assumptions there state an additional requirement that the derivative is uniformly bounded away from the origin. This requirement is not needed to obtain the generalization in this paper.

The following assumption strengthens the requirements on the transition dynamics~\eqref{eq:filters:transitions} and weakens the requirements on the observation model~\eqref{eq:filters:observations}.
\begin{assumption} \label{ass:filters:assumption_2}
    The following statements hold:
    \begin{enumerate}[label=\alph*)]
        \item[\rm{(i)}] \label{item:filters:F2} $\Xb=\Rb^d,$ $\Xc$ is an open subset of $\Rb^d,$ $\mu \ll \lambda^{[d]},$ and the function $((x,a),\xi) \mapsto F(x,a,\xi)$ satisfies the Diffeomorphic Condition;
        \item[\rm(ii)] \label{item:filters:G2} either $G$ does not depend on $a,$ or the following conditions {\rm(a)} and {\rm(b)} hold: {\rm(a)} $\Yb=\Rb^m,$ $\Hc$ is an open subset of $\Rb^m,$ and $\nu \ll \lambda^{[m]},$ and {\rm(b)} for each $x \in \Xb$ the function $(a,\eta) \mapsto G(a,x,\eta)$ satisfies the Diffeomorphic Condition.
    \end{enumerate}
\end{assumption}

As shown in Section~\ref{sec:kernels}, below, Assumption~\ref{ass:filters:assumption_2}(i) 
implies continuity in total variation of the transition kernel $\Tc,$ and Assumption~\ref{ass:filters:assumption_2}(ii) 
implies continuity in total variation in $a$ of the observation kernel $Q.$
The following theorem, whose proof follows from Theorems~\ref{thm:feller:semi-uniform_Feller_weak_continuous_filter} and \ref{thm:feller:filter_kernel_continuous}, is the main result of this paper. In particular, it generalizes the results of~\citet[Section 8.1]{feinberg_partially_2016} to models with finite dimensional state, action, observation, and noise spaces.
\begin{theorem}[{Weak continuity of the filter}] \label{thm:filters:main_result}
    Under Assumption~\ref{ass:filters:assumption_1} or ~Assumption~\ref{ass:filters:assumption_2}, the stochastic kernel $q$ on $\Pb(\Xb)$ given $\Pb(\Xb) \times \Ab$ defined in~\eqref{eq:filters:filter_kernel} is weakly continuous.
\end{theorem}

\begin{proof}
See~Section~\ref{sec:proofs}.
\end{proof}

\section{Semi-Uniform Feller Kernels and Their Properties} \label{sec:feller}

In this section we consider continuity properties of the stochastic kernel $P$ defined in~\eqref{eq:filters:joint_kernel}. The form of continuity of interest to us is semi-uniform Feller continuity, introduced by~\citet{feinberg_markov_2022,feinberg_2023_semi}. The main result is Theorem~\ref{thm:feller:semi-uniform_Feller_weak_continuous_filter}, which shows that $P$ is semi-uniform Feller if either Assumptions \ref{ass:filters:assumption_1} or \ref{ass:filters:assumption_2} hold. Suppose $\Psi$ is a stochastic kernel on $\Sb_1 \times \Sb_2$ given $\Sb_3.$ For $A \in \Bc(\Sb_1),$ $B \in \Bc(\Sb_2),$ and $s_3 \in \Sb_3,$ let
\begin{equation}\label{eq:feller:joint_marginals}
    \Psi(A,B|s_3) \defeq \Psi(A\times B|s_3).
\end{equation}
In particular, we consider marginal stochastic kernels
$\Psi(\Sb_1,\:\cdot\:|\:\cdot\:)$ on $\Sb_2$ given $\Sb_3$ and $\Psi(\:\cdot\:,\Sb_2|\:\cdot\:)$ on $\Sb_1$ given~$\Sb_3.$

\begin{definition}[{\citet[Definition 4.1]{feinberg_markov_2022},~\cite[Definition 1]{feinberg_2023_semi}}] \label{def:feller:semi-uniform_feller}
    A stochastic kernel $\Psi$ on $\Sb_1\times\Sb_2$ given $\Sb_3$ is \emph{semi-uniform Feller} if, for each $s_3 \in \Sb_3$ and for each bounded continuous function $f:\Sb_1\to\Rb,$
    \begin{equation}\label{eq:feller:semi-uniform_feller}
        \lim_{s_3' \to s_3}
        \sup_{B\in \Bc(\Sb_2)} \abs*{
            \int_{\Sb_1} f(s_1) \ \Psi(ds_1,B|s_3')-\int_{\Sb_1} f(s_1) \ \Psi(ds_1,B|s_3)
        } = 0.
    \end{equation}
\end{definition}
A semi-uniform Feller stochastic kernel $\Psi$ on $\Sb_1\times \Sb_2$ given $\Sb_3$ is weakly continuous \cite[Lemma 4.2]{feinberg_markov_2022}.

The remainder of this section describes sufficient conditions for semi-uniform Feller continuity of $P.$ According to the following statement, this property implies weak continuity of the stochastic kernel $q$ in~\eqref{eq:filters:filter_kernel}.

\begin{theorem}[{\citet[Corollary 6.7]{feinberg_markov_2022}}] \label{thm:feller:semi-uniform_Feller_weak_continuous_filter}
    If the stochastic kernel $P$ on $\Xb \times \Yb$ given $\Xb \times \Ab$ defined in~\eqref{eq:filters:joint_kernel} is semi-uniform Feller, then the stochastic kernel $q$ on $\Pb(\Xb)$ given $\Pb(\Xb) \times \Ab$ defined in~\eqref{eq:filters:filter_kernel} is weakly continuous.
\end{theorem}

The next statement provides sufficient conditions for the kernels $\Tc$ and $Q$ to imply that $P$ is semi-uniform Feller.

\begin{theorem}[{\citet[Corollary 6.11]{feinberg_markov_2022}}] \label{thm:feller:semi-uniform_feller}
   For the transition and observation stochastic kernels $\Tc$ on $\Xb$ given $\Xb \times \Ab$ and $Q$ on $\Yb$ given $\Ab \times \Xb$ defined in~\eqref{eq:filters:model_kernels} and for the stochastic kernel $P$ on $\Xb \times \Yb$ given $\Xb \times \Ab$ defined in~\eqref{eq:filters:joint_kernel}, if either of the following assumptions hold:
    \begin{enumerate}
        \item[\rm{(i)}] $\Tc$ is weakly continuous, and $Q$ is continuous in total variation; \label{item:feller:transitions_weak_observations_total_variation}
        \item[\rm{(ii)}] $\Tc$ is continuous in total variation, and $Q$ is continuous in $a$ in total variation; \label{item:feller:transitions_total_variation_observations_total_variation_a}
    \end{enumerate}
    then $P$ is semi-uniform Feller.
\end{theorem}

The major significance of Theorem~\ref{thm:feller:semi-uniform_feller} for this paper is that its assumption~(i) corresponds to Assumption~\ref{ass:filters:assumption_1}, while assumption~(ii) corresponds to Assumption~\ref{ass:filters:assumption_2}.
Therefore, Assumptions~\ref{ass:filters:assumption_1} and~\ref{ass:filters:assumption_2} 
each imply that the stochastic kernel $P$ is semi-uniform Feller. These observations are formulated in the next theorem, which together with Theorem~\ref{thm:feller:semi-uniform_Feller_weak_continuous_filter} implies Theorem~\ref{thm:filters:main_result}. 

\begin{theorem}[Semi-uniform Feller of the stochastic kernel $P$] \label{thm:feller:filter_kernel_continuous}
    Under Assumption~\ref{ass:filters:assumption_1} or~Assumption~\ref{ass:filters:assumption_2}, the stochastic kernel $P$ on $\Xb \times \Yb$ given $\Xb \times \Ab$ defined in~\eqref{eq:filters:joint_kernel} is semi-uniform  Feller. 
\end{theorem}

\begin{proof}
    See~Section~\ref{sec:proofs}.
\end{proof}

\section{Continuity of Stochastic Kernels Defined by Stochastic Equations} \label{sec:kernels}

Let $\Sb_1,$ $\Sb_2,$ and $\Omega$ again denote Borel spaces, and let us consider the stochastic kernel $\kappa$ on $\Sb_1$ given $\Sb_2.$ This section deals with two distinct but related issues concerning $\kappa.$ We shall first consider the existence of a Borel measurable function $\phi : \Sb_2 \times \Omega \to \Sb_1$ and probability measure $p \in \Pb(\Omega)$ such that the representation~\eqref{eq:notation:stochastic_kernel} holds. Theorem~\ref{thm:kernels:aumann} is Aumann's lemma~\cite[Lemma F]{aumann_1964_mixed}, which shows that such a representation exists with $\Omega = [0,1]$ and $p = \lambda^{[1]}.$ Corollary~\ref{cor:kernels:aumann_multidim} states that for every natural number $n$ such a representation exists with $\Omega = [0,1]^n$ and $p = \lambda^{[n]}.$

We will then consider conditions on $\phi$ and $p$ that imply different continuity properties of the stochastic kernel $\kappa.$ Theorem~\ref{thm:kernels:weak} shows that a necessary and sufficient condition for $\kappa$ to be weakly continuous is for $\phi$ to be continuous in distribution $p.$ It also shows that stochastic continuity is a sufficient condition for $\phi$ to imply that $\kappa$ is weakly continuous. Theorem~\ref{thm:kernels:total_variation} gives a necessary and sufficient condition for $\kappa$ to be continuous in total variation; namely, that $\phi$ is continuous in total variation {\it wrt} $p.$ In addition, it shows that if $p$ is a continuous distribution on $\Rb^n$ and $\phi$ satisfies the Diffeomorphic Condition, then $\kappa$ is continuous in total variation. Finally, we provide examples showing the significance of the assumption that $p$ is a continuous distribution for continuity in total variation.

The significance of these representation and continuity theorems for the analysis of POMDPs in this paper is that they provide the justification for Assumptions~\ref{ass:filters:assumption_1} and~\ref{ass:filters:assumption_2} each to imply one of the hypotheses of Theorem~\ref{thm:filters:main_result}. As a consequence, either assumption implies weak continuity of the filter kernel $q$ on $\Pb(\Xb)$ given $\Pb(\Xb) \times \Ab$ defined in~\eqref{eq:filters:filter_kernel}.

The following theorem is due to~\citet[Lemma F]{aumann_1964_mixed}; see also~\citet[Lemma 1.2]{gihman_1979_controlled}. It shows that $\kappa$ always admits a representation~\eqref{eq:notation:stochastic_kernel} when $\Omega = [0,1]$ and $p$ is the uniform distribution.

\begin{theorem}[{\citet[Lemma F]{aumann_1964_mixed},~\citet[Lemma 1.2]{gihman_1979_controlled}}] \label{thm:kernels:aumann}
    Let $\Sb_1$ and $\Sb_2$ be Borel spaces, and let $\kappa$ be a stochastic kernel on $\Sb_1$ given $\Sb_2.$  Then there exists a Borel measurable function  $\phi : \Sb_2 \times [0,1] \to \Sb_1,$ where the Borel $\sigma$-algebra is considered on the unit interval $[0,1],$ such that
    \begin{align} \label{eq:kernels:kernel_representation}
        \kappa(B|s_2) &= \int_0^1 \mathbf{1}\{\phi(s_2,\omega) \in B\} \ d\omega,
        &&
        B \in \Bc(\Sb_1).
    \end{align}
\end{theorem}


As a remark, we would like to mention a relevant fact, which is not used in this paper. {
\citet[Theorem 1.1]{kifer_1986_ergodic} and the isomorphism theorem for Borel spaces imply that it is possible to define a probability measure $\mathfrak{m}$ on the space $\mathfrak{B}(\Sb_2, \Sb_1)$ of Borel measurable functions $f:\Sb_2 \to \Sb_1$ such that $\kappa(B|s_2) = \int_{\mathfrak{B}(\Sb_2,\Sb_1)} \mathbf{1}\{f(s_2) \in B\} \ \mathfrak{m}(df).$}

The following {corollary} is a generalization of Theorem~\ref{thm:kernels:aumann}. It shows that the space $[0,1]$ can be replaced by $[0,1]^n,$ { $n=1,2,\ldots,$ }  without loss of generality.

\begin{corollary} \label{cor:kernels:aumann_multidim}
    Let $\Sb_1$ and $\Sb_2$ be Borel spaces, and let $\kappa$ be a stochastic kernel on $\Sb_1$ given $\Sb_2.$ Then for each natural number $n$ there exists a Borel measurable function $\phi : \Sb_2 \times [0,1]^n \to  \Sb_1,$ where the Borel $\sigma$-algebra is considered on the unit box $[0,1]^n,$ such that 
    \begin{align} \label{eq:kernels:multidim}
        \kappa(B|s_2) &= \int_{[0,1]^n} \mathbf{1}\{\phi(s_2, \omega) \in B\} \ d\omega,
        && B \in \Bc(\Sb_1).
    \end{align} 
\end{corollary}

\begin{proof}
    See Section~\ref{sec:proofs}.
\end{proof}

The remainder of this section concerns continuity of the stochastic kernel $\kappa$ defined in \eqref{eq:notation:stochastic_kernel}. There are well-known sufficient conditions on the function $\phi$ that imply that $\kappa$ is weakly continuous. For example, {continuity of $\phi$ is sufficient; \citet[p. 92]{hernandez-lerma_adaptive_1989}.} As discussed in~\citet[Section 8.1]{feinberg_partially_2016}, it is  sufficient for $s_2 \mapsto \phi(s_2, \omega)$ to be continuous for $p$-a.s. $\omega.$ The following theorem, whose proof follows directly from the corresponding definitions of continuity, describes necessary and sufficient conditions on the Borel measurable function $\phi$ for the kernel $\kappa$ to be weakly continuous. It also shows that stochastic continuity is another sufficient condition for weak continuity of $\kappa,$ which is weaker than $p$-a.s. continuity of $\phi$ at each $s_2 \in \Sb_2,$ which in turn is weaker than continuity of $\phi.$

\begin{theorem}[Weak continuity] \label{thm:kernels:weak}
    Let $p \in \Pb(\Omega)$ and $\phi : \Sb_2 \times \Omega \to \Sb_1$ be Borel measurable, and consider the stochastic kernel $\kappa$ on $\Sb_1$ given $\Sb_2$ defined in~\eqref{eq:notation:stochastic_kernel}. The following statements hold:
    \begin{enumerate}
        \item[\rm(a)] the function $\phi$ is continuous in distribution $p$ if and only if the stochastic kernel $\kappa$ is weakly continuous;
        \item[\rm(b)] if the function $\phi$ is continuous in probability $p,$ then $\kappa$ is weakly continuous.
    \end{enumerate}
\end{theorem}

\begin{proof}
    The proof is found in Section~\ref{sec:proofs}.
\end{proof}

The next theorem states necessary and sufficient conditions for the stochastic kernel $\kappa$ to be continuous in total variation. It also shows that  for  $\Sb_1 = \Rb^n$ and for $\Omega$ being an open subset of $\Rb^n,$  if $p$ is a continuous distribution and $\phi$ satisfies the Diffeomorphic Condition, the stochastic kernel $\kappa$ is continuous in total variation.  While Theorem~\ref{thm:kernels:total_variation}(a) follows directly from definitions of convergence in total variation, the proof of Theorem~\ref{thm:kernels:total_variation}(b) is nontrivial, and it uses Lemma~\ref{lem:proofs:images_converge_in_measure}. For $p \in \Pb(\Omega),$ where $\Omega$ is a Borel subset of $\Rb^n,$ we write $p\ll \lambda^{[n]}$ if $p$ is absolutely continuous {\it wrt}  $\lambda^{[n]}$ restricted to $\Omega.$

\begin{theorem}[Continuity in total variation] \label{thm:kernels:total_variation}
    Let $p \in \Pb(\Omega)$ and $\phi : \Sb_2 \times \Omega \to \Sb_1$ be Borel measurable, and consider the stochastic kernel $\kappa$ on $\Sb_1$ given $\Sb_2$ defined in~\eqref{eq:notation:stochastic_kernel}. The following statements hold:
    \begin{enumerate}
        \item[\rm(a)] \label{item:kernels:total_variation:iff} the function $\phi$ is continuous in total variation {\it wrt} $p$ if and only if the stochastic kernel $\kappa$ is continuous in total variation;
        \item[\rm(b)] \label{item:kernels:total_variation:differentiability_condition} if $\Sb_1 = \Rb^n,$ $\Omega$  is an open subset of $\Rb^n,$ $p \ll \lambda^{[n]},$ 
            and the function $\phi$ satisfies the Diffeomorphic Condition, then the stochastic kernel $\kappa$ is continuous in total variation.
    \end{enumerate}
\end{theorem}

\begin{proof}
    The proof of Theorem~\ref{thm:kernels:total_variation} can be found in Section~\ref{sec:proofs}.
\end{proof}

The next two counterexamples demonstrate the importance of the assumption $p \ll \lambda^{[n]}$ in Theorem~\ref{thm:kernels:total_variation}.


\begin{example}[Discrete distributions failing continuity in total variation] \label{ex:kernels:delta_measure}
    Suppose $\Sb_1=\Sb_2=\Omega=\Rb,$ suppose the distribution $p$ is concentrated at the point $0,$ that is, $p(B)=\mathbf{1}\{0 \in B\}$ for $B \in \Bc(\Rb),$  and suppose $\phi(s_2,\omega)=s_2+\omega.$  Then $\kappa(B|s_2) = \mathbf{1}\{s_2 \in B\}$ for $B \in \Bc(\Rb)$ since $p(0)=1.$ Consider the sequence $s_2^{(k)} = k^{-1}$ for $k=1,2,\dots,$ and let $B = \set{0}.$ Then $s_2^{(k)} \to 0,$ $\kappa(B|0) = 1,$ and $\kappa(B|s_2^{(k)}) = 0$ for all $k=1,2,\dots;$ hence, $\kappa$ is not continuous in total variation.
\end{example}


\begin{example}[Singular distributions failing continuity in total variation] \label{ex:kernels:degenerate_gaussian}
    Let $\Sb_1=\Rb^n,$ $\Sb_2=\Rb,$  $\Omega = \Rb^n,$ where $n=2,3,\ldots,$ and consider a white noise process $p = N(0, \Sigma),$ where the covariance matrix $\Sigma \in \Rb^{n \times n}$ is assumed to be singular. Then $p$ is supported on a proper subspace $V$ of $\Omega.$ Let $z \in V^\perp$ be a nonzero vector, and let us consider the kernel $\kappa$ on $\Rb^n$ given $\Rb$ defined in
    \begin{align*}
        \kappa(B|s_2) &=
        \int_{\Omega} \mathbf{1}\{s_2 z + \omega \in B\} \ p(d\omega), && B \in \Bc(\Rb^n), \qquad s_2 \in \Rb.
    \end{align*}
    We observe that $\kappa(V|s_2) = 0$ if $s_2 \ne 0,$ but $\kappa(V|s_2) = 1$ if $s_2 = 0.$ Thus, $\kappa$ is not continuous in total variation.
\end{example}


Example~\ref{ex:kernels:degenerate_gaussian} shows that  embedding a Gaussian distribution in a higher-dimensional space can cause failure of continuity in total variation.  {An} important implication of Example~\ref{ex:kernels:degenerate_gaussian}  for additive Gaussian noise models {considered in Section~\ref{sec:applications} is that} continuity in total variation of the observation kernel $Q$ can fail when the observation noise distribution is degenerate. {Thus, in general,} noise cannot be shared across vector components.

\section{Existence of Optimal Policies, Validity of Optimality Equations, and Convergence of Value Iterations for POMDPs} \label{sec:pomdps}
In view of Aumann's lemma (Theorem~\ref{thm:kernels:aumann}), there are two equivalent ways to define a POMDP: (ii) by using transition and observation kernels and (ii) by using stochastic equations \eqref{eq:filters:model}.  In this section we consider an application of the results of Sections~\ref{sec:filters} and~\ref{sec:feller} and the results on POMDPs from \cite{feinberg_markov_2022} for MDPs defined by transition and observation kernels to derive sufficient conditions for the existence of optimal policies, validity of optimality equations, and convergence of value iterations for total-cost POMDPs defined by stochastic equations \eqref{eq:filters:model}.

In this section we consider a POMDP defined by a tuple $(\Xb, \Yb, \Ab, \Tc, Q, c),$ where $\Xb,$ $\Yb,$ and $\Ab$ are the Borel spaces consisting of states, observations, and actions, $\Tc$ is the transition kernel on $\Xb$ given $\Xb \times \Ab$ defined in~\eqref{eq:filters:transition_kernel}, $Q$ is the observation kernel on $\Yb$ given
$\Ab \times \Xb$ defined in~\eqref{eq:filters:observation_kernel}, and $c : \Xb \times \Ab \to (-\infty, +\infty]$ is a bounded-below Borel measurable one-step cost function. The POMDP evolves as follows. The initial unobservable state $x_0$ is assumed to have the distribution $p_0 \in \Pb(\Xb).$ The initial observation follows
$y_0 \sim Q_{0}\left(\:\cdot\:| x_{0}\right),$ where $Q_0$ is the initial observation kernel on $\Yb$ given $\Xb.$ At each time epoch $t=0,1, \ldots,$ if the state of the system is $x_{t} \in \Xb,$ the decision maker observes $y_t$ and chooses an action $a_{t} \in \Ab,$ where $y_{t} \sim Q(\:\cdot\: | a_{t-1}, x_{t})$ for $t > 0,$ the cost $c(x_t, a_t)$ is incurred, and the system moves to the state $x_{t+1} \sim \Tc(\:\cdot\:|x_t,a_t).$

Let $\Hb_0 \defeq \Pb(\Xb) \times \Yb$ and $\Hb_{t+1} \defeq \Hb_t \times (\Ab \times \Yb)$ denote the set of observable histories at times $t=0,1,\dots.$ A policy $\pi$ for the POMDP is a sequence $\set{\pi_t}_{t=0}^{\infty}$ of stochastic kernels $\pi_t$ on $\Ab$ given $\Hb_t.$ Let $\Pi$ denote the set of policies. A policy $\pi$ is nonrandomized if each kernel $\pi_t(\:\cdot\:|h_t)$ ($h_t \in \Hb_t$) is concentrated on a single point; $\pi$ is Markov if each $\pi_t$ only depends on the current state and time; $\pi$ is stationary if each $\pi_t$ only depends on the current state. Let $\Hb \defeq (\Xb \times \Yb \times \Ab)^{\infty}$ denote the set of trajectories endowed with the product $\sigma$-algebra. The Ionescu Tulcea theorem (see, e.g.,~\citet[Proposition 7.28]{bertsekas_stochastic_1996}) implies that for every policy $\pi,$ initial distribution
$p_0 \in \Pb(\Xb),$ and stochastic kernels $\Tc,$ $Q,$ and $Q_0,$ there exists a unique probability measure $\mathbf{P}_{p_0}^\pi$ on $\Hb.$ Expectation  {\it wrt} this measure is denoted by $\mathbf{E}_{p_0}^{\pi}.$ For a horizon $T=1,2,\dots, \infty$ and discount rate $\alpha \geq 0,$ the expected total discounted costs are
\begin{align} \label{eq:pomdps:value_fn}
    V_{T,\alpha}^{\pi}(p_0) &\defeq \mathbf{E}_{p_0}^\pi \sum_{t=0}^{T-1} \alpha^t c(x_t, a_t), &&
    p_0 \in \Pb(\Xb), \quad \pi \in \Pi,
\end{align}
where $V_{0,\alpha}^{\pi}(p_0) \defeq 0.$ We consider two assumptions on the discount $\alpha$ and the cost function $c$ that each imply that~\eqref{eq:pomdps:value_fn} is well-defined.

\noindent{\bf Assumption (\textbf{D})}
    \emph{$c$ is bounded below on $\Xb \times \Ab$ and $\alpha \in [0,1).$}

\noindent{\bf Assumption (\textbf{P})}
    \emph{$c$ is nonnegative on $\Xb \times \Ab$ and $\alpha\ge 0 .$}

For Assumption~(\textbf{D}) there is no loss in generality to assuming that $c$ is nonnegative. An MDP can be viewed as a special case of a POMDP when $\Yb = \Xb$ and $Q(x_{t+1}|a_t, x_{t+1}) = 1.$ For an MDP, we shall denote the objective function in  ~\eqref{eq:pomdps:value_fn}   by $v_{T,\alpha}^{\pi}$ and the corresponding value function by $v_{T,\alpha}$ in lowercase to distinguish them from notations for  POMDPs.

A POMDP can be reduced to a \emph{completely-observable MDP} (\comdp{}) whose states in the reduced process are probability measures on the state space of the POMDP. In particular, the POMDP $(\Xb, \Yb, \Ab, \Tc, Q, c)$ can be reduced to the \comdp{} $(\Pb(\Xb), \Ab, q, \bar c),$ where the transition kernel $q$ on $\Pb(\Xb)$ given $\Pb(\Xb) \times \Ab$ is defined in~\eqref{eq:filters:filter_kernel}, and the one-step cost function $\bar{c}: \Pb(\Xb) \times \Ab \to (-\infty, +\infty]$ given by the equation
\begin{align} \label{eq:pomdps:filter_cost}
    \overline{c}(z,a) &\defeq
    \int_{\Xb}c(x,a) \ z(dx), &&
    z \in \Pb(\Xb), \quad
    a \in \Ab,
\end{align}
and $\bar c$ is well-defined since $c$ is bounded below and Borel measurable. Through a similar construction one arrives at objective criteria of the form~\eqref{eq:pomdps:value_fn} for the \comdp{} $(\Pb(\Xb), \Ab, q, \bar c).$ For the \comdp{}, histories are $(z_0, a_0, \dots, z_{t-1}, a_{t-1}, z_t)$ and trajectories are $(z_0, a_0, z_1, a_1, \dots)$ for measures $z_j \in \Pb(\Xb)$ ($j=0,1,\dots$), where $z_0=p_0.$ For more details on this reduction, see~\citet[Section 10.3]{bertsekas_stochastic_1996},~\citet[Chapter 8]{dynkin_1979_controlled},~\citet[Chapter 4]{hernandez-lerma_adaptive_1989},~\citet{rhenius_1974}, or~\citet{yushkevich_1976_reduction}. Denote $\bar{v}_{T,\alpha}$ the associated value function for $T=0,1,\dots, \infty$ for the \comdp{}.

The following assumption was introduced in~\citet{feinberg_average_cost_2012} and pertains to basic structural and continuity results for the \comdp{} $(\Pb(\Xb), \Ab, \bar c, q).$
\begin{unAssumption}[Assumption (\textbf{W*})]
    The following statements hold:
    \begin{enumerate}
        \item[\rm(i)] The function $\bar c$ is $\Kb$-inf-compact and bounded below on $\Pb(\Xb) \times \Ab.$
        \item[\rm(ii)] The stochastic kernel $q$ on $\Pb(\Xb)$ given $\Pb(\Xb) \times \Ab$ is weakly continuous.
    \end{enumerate}
\end{unAssumption}
If $c$ is $\Kb$-inf-compact on $\Xb \times \Ab$ and bounded below, then $\bar c$ is $\Kb$-inf-compact on $\Pb(\Xb);$ see~\citet[Theorem 3.4]{feinberg_partially_2016}. According to~\citet[Theorem 3.1]{feinberg_partially_2016}, if Assumption (\textbf{W*}) is
combined with either Assumption (\textbf{P}) or (\textbf{D}), then optimality equations hold for the \comdp{} $(\Pb(\Xb), \Ab, q, \bar c),$ optimal policies exist, value iterations converge, 
and the value functions $\bar v_{T,\alpha}$ are lower semicontinuous.

We are now ready to apply Theorem~\ref{thm:filters:main_result} to obtain continuity of the kernel $q$ for POMDPs with transitions and observations determined by the stochastic equations~\eqref{eq:filters:model}.  We adopt the notation
    \begin{equation}\label{eq:pomdps:integrand}
        w_{T,\alpha}(z,a) \defeq \int_{\Xb} \left[c(x,a) + \alpha
        \int_{\Xc}
        \int_{\Hc}
        \bar v_{T,\alpha} (H(z,a,G(a,F(x,a,\xi), \eta))) \
        \nu(d\eta) \ \mu(d\xi) \ \right]
        z(dx),
    \end{equation}
    for  $z \in \Pb(\Xb),$ $a \in \Ab,$ and $T = 0, 1, \dots$ or $T=\infty.$

\begin{theorem} \label{thm:pomdps:pomdps_from_equations}
    Consider the POMDP with transition and observation kernels $\Tc$ on $\Xb$ given $\Xb \times \Ab$ and $Q$ on $\Yb$ given $\Ab \times \Xb$ defined in~\eqref{eq:filters:model_kernels} for the model~\eqref{eq:filters:model}. Under either Assumption~\ref{ass:filters:assumption_1} or Assumption~\ref{ass:filters:assumption_2}, the  transition kernel $q$ on $\Pb(\Xb)$ given $\Pb(\Xb) \times \Ab$ for the \comdp{} defined in~\eqref{eq:filters:filter_kernel} is weakly continuous.  Therefore, if either Assumption (\textbf{D}) or (\textbf{P}) holds, and the one-step cost function $c$ is $\Kb$-inf-compact, then the conclusions of \cite[Theorem 3.1]{feinberg_partially_2016} hold. In particular, for each $T = 0, 1, \dots, \infty,$ the value functions $\bar{v}_{T+1,\alpha}(z)$ are lower semicontinuous, and they satisfy the optimality equations
    \begin{equation}
        \bar{v}_{T+1,\alpha}(z) =  \min_{a \in \Ab} \set*{
            w_{T,\alpha}(z,a)
        },
    \end{equation}
    where $\psi_{T,\alpha}(z,a)$ is defined in~\eqref{eq:pomdps:integrand}, and the nonempty sets of optimal actions $A_{T,\alpha}$ are given by
    \begin{equation}
        A_{T,\alpha}(z) \defeq
        \set*{a \in \Ab : \bar{v}_{T+1,\alpha}(z) =
           w_{T,\alpha}(z,a)
        }.
    \end{equation}
   { For the \comdp{} defined in~\eqref{eq:filters:filter_kernel}, there  exist nonrandomized Markov optimal policies for $T$-horizon problems}, and, if $\varphi=(\varphi_0,\varphi_1,\ldots,\varphi_{T-1})$ is a nonrandomized Markov policy for the \comdp{} such that $\varphi_t(z)\in A_{T-t-1,\alpha}(z)$ for all $t=0,1,\ldots,T-1$ and for all $z\in \Pb(\Xb),$ then the policy $\varphi$ is optimal.   For an infinite-horizon problem,
    \[
    \bar{v}_{\infty,\alpha}(z)=\lim_{T\to\infty} \bar{v}_{T,\alpha}(z), \qquad z\in \Pb(\Xb),
    \]
    and  { for the \comdp{}}   there exists a nonrandomized stationary optimal policy, and a nonrandomized stationary policy $\varphi$ { for the \comdp{}} is optimal if an only if $\varphi(z)\in A_{\infty,\alpha}(z)$ for all $z\in \Pb(\Xb).$
\end{theorem}

\begin{proof}
    This theorem follows directly from Theorems~\ref{thm:feller:semi-uniform_Feller_weak_continuous_filter}, \ref{thm:feller:semi-uniform_feller}  and \citet[Corollary 6.12]{feinberg_markov_2022}.
\end{proof}

By applying standard procedures, nonrandomized Markov optimal policies and nonrandomized  stationary optimal policies for a COMDP can be used to construct nonrandomized optimal policies for the original POMDP; see e.g., \cite[p. 90]{hernandez-lerma_adaptive_1989}.

\section{Examples of Applications} \label{sec:applications}

This section deals with the version of equations~\eqref{eq:filters:model} in which the function $G$ does not depend on the action:
\begin{subequations}\label{eq:applications:general_kalman_filter}
    \begin{align}
        x_{t+1} &= F(x_t, a_t, \xi_t), &&
        x_t \in \Xb, \quad
        a_t \in \Ab, \quad
        \xi_t \in \Xc,
        \label{eq:applications:general_kalman_filter_states} \\
        y_{t+1} &= \tilde{G}(x_{t+1}, \eta_{t+1}), &&
        y_{t+1} \in \Yb, \quad
        \eta_{t+1} \in \Hc.
        \label{eq:applications:general_kalman_filter_obs}
    \end{align}
\end{subequations}
Such equations commonly appear in statistical filtering theory, as discussed below. Model~\eqref{eq:applications:general_kalman_filter} is a simpler model than~\eqref{eq:filters:model} because it does not include actions in the observation equation.  All statements in this section deal with model \eqref{eq:applications:general_kalman_filter} with Borel measurable functions $F$ and $\tilde{G},$ and with the sets and variables defined for model
\eqref{eq:filters:model}. Recall that noises are independent, and $\xi_t\sim\mu,$  $\eta_t~\sim\nu.$ The following theorem is the main result of this section. 

\begin{theorem}[Nonlinear filtering] \label{thm:applications:state_space_models}
    For model~\eqref{eq:applications:general_kalman_filter},
   under either of the  following two assumptions (i) or (ii):
    \begin{enumerate}
        \item[\rm(i)] \label{item:applications:ssm_case_1} the following statements hold:
        \begin{enumerate}
            \item \label{item:applications:ssm_case_1a} the function $((x,a), \xi) \mapsto F(x,a,\xi)$ is continuous in distribution $\mu;$
            \item \label{item:applications:ssm_case_1b} $\Yb = \Rb^m,$ $\Hc$ is an open  subset of $\Rb^m,$ $\nu \ll \lambda^{[m]},$ and the function $(x,\eta) \mapsto \tilde{G}(x, \eta)$ satisfies the Diffeomorphic Condition;
        \end{enumerate}
        \item[\rm(ii)] \label{item:applications:ssm_case_2} $\Xb= \Rb^d,$ $\Xc$ is an open  subset of $\Rb^d,$ $\mu \ll \lambda^{[d]},$ and the function $((x,a),\xi) \mapsto F(x,a,\xi)$ satisfies the Diffeomorphic Condition;
    \end{enumerate}
 the transition probability for the filter, which is the stochastic kernel $q$ on $\Pb(\Xb)$ given $\Pb(\Xb) \times \Ab$ defined in~\eqref{eq:filters:filter_kernel}, is weakly continuous.  In addition,  if either Assumption (\textbf{D}) or (\textbf{P}) from Section~\ref{sec:pomdps}  is satisfied, and the one-step cost function $c$ is $\Kb$-inf-compact, then the conclusions of Theorem~\ref{thm:pomdps:pomdps_from_equations} hold. In particular, for  expected total costs, value functions are lower semicontinuous, finite-horizon values converge to the infinite-horizon value,  optimality equations hold, and, { for the \comdp{} described in~\eqref{eq:filters:filter_kernel}}, optimality equations define  Markov optimal policies for finite-horizon problems and stationary optimal policies for infinite-horizon problems. 
\end{theorem}

\begin{proof}
    This follows directly from Theorems~\ref{thm:filters:main_result} and~\ref{thm:pomdps:pomdps_from_equations} applied to  model~\eqref{eq:applications:general_kalman_filter}.
\end{proof}

\begin{example}\label{exa:diff}
{\rm Let us consider  problem \eqref{eq:applications:general_kalman_filter} with $F(x,a,\xi)=(x^2+1)(a^2+1)\arctan(\xi),$  $x,a,\xi\in \Rb.$ In addition, the observation function ${\tilde G}$ is  Borel measurable, and the one-step function $c$ is $\Kb$-inf-compact. All noises are Gaussian.  This example satisfies assumptions of Theorem~\ref{thm:applications:state_space_models}.  In particular, assumption (ii) holds.  Thus, the conclusions of Theorem~\ref{thm:applications:state_space_models} hold.  In addition, the images sets $F(x,a,\Rb)$ are not constant.  For example, $F(0,0,\Rb) = (-\frac{\pi}{2},\frac{\pi}{2})\ne (-\pi,\pi)=F(1,0,\Rb).$
}
\end{example}

 The following corollary covers the cases of additive and multiplicative transition and observation noises broadly used in the statistical filtering literature.
\begin{corollary}\label{corapplic} The conclusions of Theorem~\ref{thm:applications:state_space_models} hold for model~\eqref{eq:applications:general_kalman_filter} in the following cases:

(a) {\rm (additive transition noise)}  $F(x_{t},a_t,\xi_{t}) = f(x_{t},a_t) + \xi_{t},$ where  $\Xb=\Xc=\Rb^d,$ $f : \Xb\times\Ab \to \Xb$ is a continuous function, and  $\mu \ll \lambda^{[d]};$

(b) {\rm (multiplicative transition noise)}   $F(x_{t},a_t,\xi_{t}) = \diag(\xi_{t}) f(x_{t},a_t),$ 
        where   $\Xb = \Xc = \Rb^d,$  $\diag(r)$ is the diagonal matrix whose diagonal entries are formed by the vector $r,$ $f : \Xb\times\Ab \to \Xb$ 
        is a continuous function such that $f_j(x_{t},a_t))\ne 0$ for all $(x_{t},a_t)\in\Xb\times\Ab$ and for all $j=1,\ldots,d,$ and  
         $\mu \ll \lambda^{[d]};$

(c) {\rm (additive observation noise)}   $\tilde{G}(x_{t+1},\eta_{t+1}) = g(x_{t+1}) + \eta_{t+1},$ where $\Yb = \Hc = \Rb^m,$       $g : \Xb \to \Yb$ is a continuous function, 
        $\nu \ll \lambda^{[m]},$  and the function $((x,a), \xi) \mapsto F(x,a,\xi)$ is continuous in distribution $\mu,$ which takes place, for example, if $F$ is continuous;

(d)  {\rm (multiplicative observation noise)}   $\tilde{G}(x_{t+1},\eta_{t+1}) = \diag(\eta_{t+1})g(x_{t+1}),$ 
    where $\Yb = \Hc = \Rb^m,$ $g : \Xb \to \Yb$ is a continuous function such that $g_i(x_{t+1})\ne 0$ for all $x_{t+1}\in\Xb$ and for all $i=1,\ldots,m,$  
    $\nu \ll \lambda^{[m]},$  and the function $((x,a), \xi) \mapsto F(x,a,\xi)$ is continuous in distribution $\mu.$
\end{corollary}
\begin{proof}
(a) Since $\mu \ll \lambda^{[d]},$ and the Jacobian $D_\xi F(x,a,\xi) = I_{d \times d}$ is the identity matrix,  this model satisfies assumption (ii) of Theorem~\ref{thm:applications:state_space_models}. (b) Since $\mu \ll \lambda^{[d]},$ and the Jacobian $D_\xi F(x,a,\xi) = \diag( f(x_{t},a_t)),$  this model satisfies assumption (ii) of Theorem~\ref{thm:applications:state_space_models}.  (c,d) In these cases, assumption (i) of Theorem~\ref{thm:applications:state_space_models} holds.
\end{proof}

\paragraph{Linear state space models (LSSMs)}
The LSSMs are the prototypical control models with a long history of study; see, e.g., Chapter 6 of~\citet{kwakernaak_1974_linear} or Chapter 5 of \citet{hendricks_2008_linear} for overviews of the subject. In such models, the spaces $\Xb,$ $\Ab,$ $\Yb,$ $\Xc$ and $\Hc$ are all Euclidean spaces, and the transition and observation equations~\eqref{eq:applications:general_kalman_filter} take the form
    \begin{equation} \label{eq:applications:lssm}
        x_{t+1} = \vec{F}_1x_t + \vec{F}_2a_t + \xi_t, \qquad 
        y_{t+1} = \vec{G} x_{t+1} + \eta_{t+1}, 
    \end{equation}
where $\vec{F}_1,$ $\vec{F}_2,$ and $\vec{G}$ are matrices of appropriate dimensions. 

\begin{corollary}[LSSMs]
        Consider the LSSM~\eqref{eq:applications:lssm} with $\Xb = \Xc = \Rb^d$ and $\Yb = \Hc = \Rb^m.$ 
        If either $\mu \ll \lambda^{[d]}$ or $\nu \ll \lambda^{[m]},$ then the conclusions of Theorem~\ref{thm:applications:state_space_models} hold.%
\end{corollary}

\begin{proof}
   This corollary follows from Corollary~\ref{corapplic}(a,c). 
\end{proof}

For control problems with LSSM dynamics or multiplicative observation noise, a variety of cost functions $c(x,a)$ employed in applications are $\Kb$-inf-compact and bounded below. In quadratic regulator problems, the controller is faced with a cost function $c(x,a) = x^T \vec{X} x + a^T \vec{A} a,$ where $\vec{X}$ and $\vec{A}$ are matrices of appropriate dimension. We recall that a matrix $\vec{M}$ is positive-definite (semidefinite) if the inequality $x^T\vec{M}x > 0$ ($x^T \vec{M} x \geq 0$) for all $x \ne 0.$ If $\vec{M}$ is positive definite, then there exists $r > 0$ such that $x^T \vec{M} x \geq r x^Tx > 0$ for all $x \ne 0,$  and it is possible to set $r := \min\set{x^T \vec{M} x : \norm{x} = 1} > 0$. If  $\vec{X}$ and $\vec{A}$ are positive definite matrices, then the function $c$ is inf-compact since $x^T \vec{X} x$ and $a^T \vec{A} a$ are inf-compact functions on $\Xb$ and $\Ab$ respectively. However, if the matrix $\vec{X}$ is positive semidefinite, and the matrix  $\vec{A}$ is positive definite, then the function $c$ may not be inf-compact because, if  $c(x,a) \leq \gamma,$ then $c(k \tilde{x},a)\le c(x,a) \leq \gamma$  for all $k\in\Rb,$  where $0\ne \tilde{x}\in\Xb$ and $\tilde{x}^T \vec{X} \tilde{x}=0,$ and the set $\{(k\tilde{x},a):k\in\Rb\},$ which is a line, cannot be a subset of a compact set.

Another common cost function arises from state estimation problems (see, e.g.,~\citet{guo_constrained_2017}), in which the controller pays the cost  $c(x,a) = \norm{\vec Xx-\vec Aa},$ where $\norm{\:\cdot\:}$ is any equivalent Euclidean norm and $\vec{A}$ is assumed to be nonsingular.  Here too, the cost function may not be inf-compact. Indeed, if there exist $\tilde{x} \ne 0$ and $\tilde{a} \ne 0$ such that $\vec X \tilde{x} = \vec A \tilde{a},$ then $\set{(x,a) : c(x,a) \leq \gamma} \supset \set{(k \tilde{x}, k \tilde{a}) : k \in \Rb},$  and the latter set, which is a line, cannot be a subset of a compact set. The following corollary shows that these cost functions are $\Kb$-inf-compact and nonnegative; hence, under the  other hypotheses of Theorem~\ref{thm:applications:state_space_models}, the conclusions of Theorem~\ref{thm:pomdps:pomdps_from_equations} hold.

\begin{corollary}[Nonlinear filtering with $\Kb$-inf-compact costs]
        Consider the spaces $\Xb = \Rb^d,$ $\Ab = \Rb^\ell,$ and suppose the cost function $c : \Xb \times \Ab \to \Rb$ is either
        \begin{enumerate}
            \item[\rm{(i)}] \label{item:applications:k_inf_compact_1} $c(x,a) \defeq x^T \vec{X}x + a^T \vec{A} a,$ where $\vec{X} \in \Rb^{d \times d}$ is positive semidefinite and $\vec{A}\in \Rb^{\ell \times \ell}$ is positive definite;
            \item[\rm{(ii)}] \label{item:applications:k_inf_compact_2} $c(x,a) \defeq \norm{\vec{X}x - \vec{A}a},$ where $\vec{X}\in \Rb^{k \times d}$ and $\vec{A}\in \Rb^{k \times \ell}$ is nonsingular.
        \end{enumerate}
        Then  the function $c$ is  nonnegative and $\Kb$-inf-compact. 
\end{corollary}

\begin{proof}
    It is obvious that $c$ is nonnegative. (i) Since $\vec{A}$ is positive definite, there exists $r > 0$ such that $a^T \vec{A} a \geq r a^Ta > 0$ for all $a \ne 0.$ Hence, if $C \subset \Xb$  is compact, then the closed set $\set{(x,a) : x \in C, c(x,a) \leq \gamma} \subset \set{(x,a) : x \in C, \norm{a} \leq \gamma r^{-1}}$ is compact, so $c(x,a)$ is $\Kb$-inf-compact. (ii) To see that $c(x,a)$ is $\Kb$-inf-compact, let $C \subset \Xb$ be compact. In particular, $C$ is contained in a closed ball of radius $r,$ so if $x \in C$ then $\norm{x} \leq r.$ Then $\norm{\vec{X}x} \leq \norm{\vec{X}} r,$ and if $c(x,a) \leq \gamma,$ then it follows that $\norm{\vec{A}a} \leq \norm{\vec{A}a - \vec{X}x} + \norm{\vec{X}x} \leq \gamma + \norm{\vec{X}} r.$ Since $\norm{\vec{A} a}$ is bounded and $\vec{A}$ is nonsingular, it follows that $\norm{a}$ is bounded; hence, the closed set $\set{(x,a) : x \in C, c(x,a) \leq \gamma}$ is compact because it is contained in a compact set, so $c(x,a)$ is indeed $\Kb$-inf-compact.
\end{proof}

\paragraph{Inventory models with noisy observations}
We consider a periodic-review inventory control problem with multiple products, continuous inventory levels, and instant orders; see, e.g.,~\cite{bensoussan2011dynamic,feinberg_2016_optimality,feinberg2022continuity} and the references therein for an overview of this problem. Let $\Xb = \Yb = \Ab = \Xc = \Hc = \Rb^d.$  At each period, the controller observes the current inventory level of each of the $d$ products, but the observation is subject to disturbances arising from measurement error~\cite{downing_1980_application}. The controller makes a determination for ordering new amounts of each product, which are fulfilled instantly. Demand for each product is realized and subtracted from the inventory level. If the problem allows for backorders, the new inventory level is simply the difference of inventory and demand. In problems with lost sales, the new inventory level for a given product does not decrease below 0. The dynamics of the system follow~\eqref{eq:applications:general_kalman_filter} with the specific functional form
    \begin{equation} \label{eq:applications:inventory_model}
        x_{t+1} = L(x_{t} + a_{t} - \xi_t), \qquad
        y_{t+1} = x_{t+1} + \eta_{t+1},
    \end{equation}
where $\set{\xi_t}_{t=0}^{\infty}$ is an iid sequence of nonnegative demands with common distribution $\mu \in \Pb(\Xb),$ $L : \Xb \to \Xb$ is defined in $L(x) = x$ for the problem with backorders and $[L(x)]^{(j)} = \max\set{0,x^{(j)}}$ (for component $j=1,\dots,d$) for problems with lost sales, and $\set{\eta_t}_{t=0}^{\infty}$ is an iid sequence of observation disturbances with common distribution $\nu \in \Pb(\Xb).$ The model~\eqref{eq:applications:inventory_model} is a special case of an additive observation noise model.

The controller faces the one-stage cost function
\begin{equation} \label{eq:applications:inventory_costs}
    c(x,a) \defeq
    \tilde{c}(a) + \int_{\Xc}  h(L(x+a-\xi)) \ \mu(d\xi)
\end{equation}
where $\tilde{c} : \Ab \to \Rb$ is an inf-compact function representing order  costs, and $h : \Xb \to [0,+\infty)$ is a convex function 
representing holding and backorder costs, which  in the case of backorders satisfies $\int h(x-\xi) \ \mu(d\xi) < +\infty$ for all $x \in \Xb.$ A typical form of the function $\tilde{c}$ is $\tilde{c}(a) \defeq \sum_{j=1}^{d} [K^{(j)}\mathbf{1}\{a^{(j)} > 0\} + \bar c^{(j)} a^{(j)}],$ where the vectors $K, \bar c \in \Xb$ have nonnegative components and represent fixed and unit order costs, respectively. The following corollary applies Theorem~\ref{thm:applications:state_space_models} to inventory control problems.

\begin{corollary}[Inventory control]
        Consider the inventory control problem~\eqref{eq:applications:inventory_model} for $\Xb = \Ab= \Yb = \Xc = \Hc = \Rb^d$ with costs~\eqref{eq:applications:inventory_costs}.
        If $\nu \ll \lambda^{[d]},$ then the conclusions of Theorems~\ref{thm:applications:state_space_models} and~\ref{thm:pomdps:pomdps_from_equations} hold.%
\end{corollary}

\begin{proof}
   This corollary follows from Corollary~\ref{corapplic}(c). The one-stage cost function $c$ is $\Kb$-inf-compact as a sum of an inf-compact function that does not depend on $x$ and a nonnegative lower semicontinuous function; see also \citet[Theorem 5.3(2)]{FL2018}.
\end{proof}


\section{%
    Another Form of Stochastic Equations%
} \label{sec:pomdp1s}

In this section we consider a variation of the model~\eqref{eq:filters:model} under different observation constraints. This leads to a different kind of partially observable controlled process, which we  call a \pomdpOne{}, {and which differs} from the POMDPs studied in the previous sections of this paper.  In \pomdpOne{} the observation $y_{t+1}$ depends directly on $(x_t,a_t),$ and it is conditionally independent on the new state $x_{t+1}$ given $(x_t,a_t).$

\citet{platzman_1980} observed that in addition to the model of POMDPs broadly used in electrical engineering and studied in the preceding sections of this paper, the \pomdpOne{} is sometimes considered in the operations research literature; see, e.g.,~\citet{smallwood_1973_optimal} and~\citet{sondik_1978_optimal}.~\citet{feinberg_markov_2022} considered both the \pomdpOne{} and the POMDP (called a \pomdpTwo{} in \cite{feinberg_markov_2022}) as special cases of a general MDP with incomplete information (MDPII); see Figure 1 in \cite{feinberg_markov_2022},  which describes the relationship between several models with incomplete information. The results for the \pomdpOne{} model are different from the results for POMDPs because, as shown in~\citet[Corollary 6.10]{feinberg_markov_2022}, weak continuity of the transition kernel and continuity of the observation kernel in total variation are the necessary and sufficient conditions for the validity of the semi-uniform Feller continuity property for \pomdpOne{}s.  However,~Theorem~\ref{thm:feller:semi-uniform_feller} provides only sufficient conditions for the validity of semi-uniform Feller continuity of the kernel $P$ for POMDPs. \citet[Example 1]{feinberg_2023_equivalent} shows that these sufficient conditions are not necessary.

We now assume that the model is given by the equations
\begin{subequations} \label{eq:pomdp1s:model}
    \begin{align}
        \label{eq:pomdp1s:transitions}
        x_{t+1} &= F(x_t, a_t, \xi_t),
        &&
        x_t \in \Xb, \quad
        a_t \in \Ab, \quad
        \xi_t \in \Xc, \\
        \label{eq:pomdp1s:observations}
        y_{t+1} &= G_1(x_t, a_t, \eta_t),
        &&
        y_{t+1} \in \Yb, \quad \eta_t \in \Hc,
    \end{align}
\end{subequations}
where the functions $F : \Xb \times \Ab \times \Xc \to \Xb$ of transition dynamics and $G_1 : \Xb \times \Ab \times \Hc \to \Yb$ of observations are Borel measurable. We observe that $F$ in~\eqref{eq:pomdp1s:transitions} agrees with~\eqref{eq:filters:transitions}, while the main difference between $G_1$ in~\eqref{eq:pomdp1s:observations} and $G$ in~\eqref{eq:filters:observations} is that the former depends on the current state $x_t,$ while the latter depends on the next state $x_{t+1}.$ Analogous to~\eqref{eq:filters:model_kernels}, we  consider the stochastic kernels
\begin{subequations}
    \label{eq:pomdp1s:model_kernels}
    \begin{align}
        \Tc(B|x, a) &= \int_{\Xc} \mathbf{1}\{F(x, a, \xi) \in B\} \ \mu(d\xi), && B \in \Bc (\Xb), \quad x \in \Xb, \quad a \in \Ab, \label{eq:pomdp1s:transition_kernel} \\
        Q_1(C|x, a) & =\int_{\Hc}\mathbf{1}\{G_1(x, a, \eta) \in C\} \ \nu(d\eta),
        && C \in \Bc (\Yb), \quad x \in \Xb, \quad a \in \Ab,  \label{eq:pomdp1s:observation_kernel}
    \end{align}
\end{subequations}
where we note that $\Tc$ in~\eqref{eq:pomdp1s:transition_kernel} agrees with~\eqref{eq:filters:transition_kernel}. 
As described in~\citet{feinberg_markov_2022}, the stochastic kernels in~\eqref{eq:pomdp1s:model_kernels} give rise to the joint stochastic kernel $P$ on $\Xb \times \Yb$ given $\Xb \times \Ab$ defined by the formula
\begin{align} \label{eq:pomdp1s:joint_kernel}
    P(B\times C|x, a) &= \Tc(B|x,a) Q_1(C|x,a), && B \in \Bc(\Xb), \quad C \in \Bc(\Yb), \quad x \in \Xb, \quad a \in \Ab.
\end{align}
Formulae~\eqref{eq:filters:joint_belief_kernel},~\eqref{eq:filters:joint_belief_conditional_representation}, and~\eqref{eq:filters:filter_kernel} applied to the stochastic kernel $P$ defined in \eqref{eq:pomdp1s:joint_kernel} implies the existence of the map $H : \Pb(\Xb) \times \Ab \times \Yb \to \Pb(\Xb)$ and the filter kernel $q$ on $\Pb(\Xb)$ given $\Pb(\Xb) \times \Ab,$ each describing the dynamics of the filter process. The following theorem, proven in~\citet{feinberg_markov_2022}, provides necessary and sufficient conditions for the stochastic kernel $P$ to be semi-uniform Feller.

\begin{theorem}[{\citet[Corollary 6.10]{feinberg_markov_2022}}] \label{thm:pomdp1s:semi-uniform_feller}
    For a \pomdpOne{} $(\Xb, \Yb, \Ab, \Tc, Q_1, c),$ the following two conditions holding together:
    \begin{enumerate}
        \item[\rm{(i)}] the stochastic kernel $\Tc$ on $\Xb$ given $\Xb \times \Ab$ defined in~\eqref{eq:pomdp1s:transition_kernel} is weakly continuous;
        \item[\rm{(ii)}] the stochastic kernel $Q_1$ on $\Yb$ given $\Xb \times \Ab$ defined in~\eqref{eq:pomdp1s:observation_kernel} is continuous in total variation;
    \end{enumerate}
    are equivalent to semi-uniform Feller continuity of the stochastic kernel $P$ on $\Xb \times \Yb$ given $\Xb \times \Ab$ defined in~\eqref{eq:pomdp1s:joint_kernel}.
\end{theorem}

In view of Theorem~\ref{thm:pomdp1s:semi-uniform_feller}, we can use Theorem~\ref{thm:filters:main_result} to give sufficient conditions for the stochastic kernel $P$ to be semi-uniform Feller continuous when $\Tc$ and $Q_1$ are generated by the model~\eqref{eq:pomdp1s:model}. These conditions imply that the kernel $q$ is weakly continuous, and are presented in the next theorem.

\begin{theorem} \label{thm:pomdp1s:results}
    If the conditions
    \begin{enumerate}
        \item[\rm{(i)}] the function $((x,a), \xi) \mapsto F(x,a,\xi)$ is continuous in distribution $\mu;$
        \item[\rm{(ii)}] $\Yb = \Rb^m,$ $\Hc$ is an open  subset of $\Rb^m,$ $\nu \ll \lambda^{[m]},$ and the function $((x,a), \eta) \to G_1(x,a,\eta)$ satisfies the Diffeomorphic Condition;
    \end{enumerate}
    hold, then the stochastic kernels $\Tc$ on $\Xb$ given $\Xb \times \Ab$ defined in~\eqref{eq:pomdp1s:transition_kernel} and $Q_1$ on $\Yb$ given $\Xb \times \Ab$ defined in~\eqref{eq:pomdp1s:observation_kernel} are weakly continuous and continuous in total variation, respectively. Hence, the stochastic kernel $P$ on $\Xb \times \Yb$ given $\Xb \times \Ab$ defined in~\eqref{eq:pomdp1s:joint_kernel} is semi-uniform Feller.
\end{theorem}

\begin{proof}
    The proof of Theorem~\ref{thm:pomdp1s:results} is provided in Section~\ref{sec:proofs}.
\end{proof}

We remark that the conditions in Theorem~\ref{thm:pomdp1s:results} simply state Assumption~\ref{ass:filters:assumption_1} for the \pomdpOne{}. Because Theorem~\ref{thm:pomdp1s:semi-uniform_feller} gives necessary and sufficient conditions for $P$ to be semi-uniform Feller, {for the \pomdpOne{} Assumption~\ref{ass:filters:assumption_2} does not imply that the stochastic kernel $P$ is semi-uniform Feller.}

    Theorem~\ref{thm:pomdp1s:results} gives sufficient conditions for the kernel $P$ to be semi-uniform Feller, which implies that the kernel $q$ is weak{ly} continuous. We close this section with a discussion of this fact in view of the results of~Section~\ref{sec:pomdps}. Let us again consider the one-step cost function $c : \Xb \times \Ab \to (-\infty,+\infty]$ and the corresponding filter cost function $\bar c : \Pb(\Xb) \times \Ab \to (-\infty, +\infty]$ defined in~\eqref{eq:pomdps:filter_cost}. If $c$ is $\Kb$-inf-compact and bounded below, then~\cite[Theorem 3.4]{feinberg_partially_2016} states that $\bar c$ is also $\Kb$-inf-compact and bounded below. Therefore, under the hypotheses of~Theorem~\ref{thm:pomdp1s:results}, the \comdp{} $(\Pb(\Xb),\Ab, q, \bar{c})$ described in~Section~\ref{sec:pomdps}, but for the transition probability $P$ defined in \eqref{eq:pomdp1s:joint_kernel}, satisfies Assumption~(\textbf{W*}). Let $\bar{v}_{T,\alpha}$ denote the value function of the \comdp{} as described in~Section~\ref{sec:pomdps} for $T = 0, 1, \dots, \infty$ and discount rate $\alpha \geq 0.$ We state an analogue of~Theorem~\ref{thm:pomdps:pomdps_from_equations} for the \pomdpOne{}.
    \begin{theorem} \label{thm:pomdp1s:pomdp1s_from_equations}
        Consider the \pomdpOne{} with transition and observation kernels $\Tc$ on $\Xb$ given $\Xb \times \Ab$ and $Q_1$ on $\Yb$ given $\Xb \times \Ab$ defined in~\eqref{eq:pomdp1s:model_kernels} for the model~\eqref{eq:pomdp1s:model}. Under the hypotheses of~Theorem~\ref{thm:pomdp1s:results}, the \comdp{} transition kernel $q$ on $\Pb(\Xb)$ given $\Pb(\Xb) \times \Ab,$ defined in~\eqref{eq:filters:filter_kernel} for the transition probability $P$ in \eqref{eq:filters:joint_belief_kernel} defined in formula \eqref{eq:pomdp1s:joint_kernel}, is weakly continuous.  Therefore, if either Assumption (\textbf{D}) or (\textbf{P}) holds, and the one-step cost function $c$ is $\Kb$-inf-compact, then the conclusions of \cite[Theorem 3.1]{feinberg_partially_2016} hold. In particular, for each $T = 0, 1, \dots, \infty,$ the value functions $\bar{v}_{T+1,\alpha}$ are lower semicontinuous, and they satisfy the optimality equations
        \begin{equation}
            \bar{v}_{T+1,\alpha}(z) =  \min_{a \in \Ab} \set*{
                \int_{\Xb} \left[c(x,a) + \alpha
                    \int_{\Hc}
                    \bar v_{T,\alpha} (H(z,a,G_1(x,a, \eta))) \
                    \nu(d\eta) \ \right]
                    z(dx)
            },
        \end{equation}
        and the nonempty sets of optimal actions $A_{T,\alpha}$ are given by
        \begin{equation}
            A_{T,\alpha}(z) \defeq
            \set*{a \in \Ab : \bar{v}_{T+1,\alpha}(z) =
                \int_{\Xb} \left[c(x,a) + \alpha
                    \int_{\Hc}
                    \bar v_{T,\alpha} (H(z,a,G_1(x,a, \eta))) \
                    \nu(d\eta) \ \right]
                    z(dx)
            }.
        \end{equation}
  For the \comdp{} defined in~\eqref{eq:filters:filter_kernel}, there  exist nonrandomized Markov optimal policies for $T$-horizon problems,
 and, if  $\varphi$ is a nonrandomized Markov policy $\varphi=(\varphi_0,\varphi_1,\ldots,\varphi_{T-1})$ for the COMDP such that $\varphi_t(z)\in A_{T-t-1,\alpha}(z)$ for all $t=0,1,\ldots,T-1$ and for all $z\in \Pb(\Xb),$ then the  policy $\varphi$ is optimal.  For an infinite-horizon problem there exists a nonrandomized stationary optimal policy  for the \comdp{}, and a nonrandomized stationary policy $\varphi$  for the \comdp{} is optimal if an only if $\varphi(z)\in A_{\infty,\alpha}(z)$ for all $z\in \Pb(\Xb).$ In addition,  $
    \bar{v}_{\infty,\alpha}(z)=\lim_{T\to\infty} \bar{v}_{T,\alpha}(z)$ for all $z\in \Pb(\Xb)$.
    \end{theorem}
    \begin{proof}
        This follows directly from Theorem~\ref{thm:pomdp1s:results}, Theorem~\ref{thm:feller:semi-uniform_Feller_weak_continuous_filter}, and
        \citet[Corollary 6.11]{feinberg_markov_2022}.
    \end{proof}

\section{Proofs of %
    Corollary~\ref{cor:kernels:aumann_multidim}
    and Theorems~\ref{thm:kernels:weak},
    \ref{thm:kernels:total_variation},
    \ref{thm:feller:filter_kernel_continuous},
    \ref{thm:filters:main_result}, and
    \ref{thm:pomdp1s:results}%
}  \label{sec:proofs}

In this technical section we provide the proofs to the main results of the paper. The outline of the proofs is  the following. Corollary~\ref{cor:kernels:aumann_multidim} and Theorems~\ref{thm:kernels:weak} and~\ref{thm:kernels:total_variation} are the principal theoretical results, and their proofs are given first. The proof of Theorem~\ref{thm:feller:filter_kernel_continuous} will follow from Theorems~\ref{thm:kernels:weak} and~\ref{thm:kernels:total_variation}, and Theorem~\ref{thm:filters:main_result} follows from Theorems~\ref{thm:feller:semi-uniform_Feller_weak_continuous_filter} and \ref{thm:feller:filter_kernel_continuous}. Finally, the proof of Theorem~\ref{thm:pomdp1s:results} follows from Theorem~\ref{thm:kernels:weak} and~\ref{thm:kernels:total_variation}.

\begin{proof}[Proof of Corollary~\ref{cor:kernels:aumann_multidim}]

    From Aumann's lemma (Theorem~\ref{thm:kernels:aumann}), there exists a function $\phi_1 : \Sb_2 \times [0,1] \to \Sb_1$ such that
    \begin{equation*}
        \kappa(B|s_2) = \int_{[0,1]} \mathbf{1}\{\phi_1(x,\omega_1) \in B\} \ d\omega_1
    \end{equation*}
    for all $B \in \Bc(\Sb_1)$ and $s_2 \in \Sb_2.$ Define $\phi : \Sb_2 \times [0,1]^n \to \Sb_1$ by the formula $\phi(s_2, \omega_1, \dots, \omega_n) \defeq \phi_1(s_2, \omega_1).$ Then
    \begin{align*}
        \int_{[0,1]^n} \mathbf{1}\{\phi(s_2, \omega_1, \dots, \omega_n) \in B\} \ d\omega_1 \cdots d\omega_n
        &=
        \int_{[0,1]^n} \mathbf{1}\{\phi_1(s_2, \omega_1) \in B\} \ d\omega_1 \cdots d\omega_n \\
        &=
        \int_{[0,1]}  \mathbf{1}\{\phi_1(s_2, \omega_1) \in B\} \ d\omega_1 = \kappa(B|s_2),
    \end{align*}
    as desired.
\end{proof}

\begin{proof}[Proof of Theorem~\ref{thm:kernels:weak}]  For statement (a), continuity of $\phi$ in distribution $p$ is equivalent to continuity of the mapping  $s_2 \mapsto \int_{\Omega} f(\phi(s_2,\omega)) \ p(d\omega),$ and  weak continuity of $\kappa$ is equivalent to continuity of the mapping  $s_2 \mapsto\int_{\Sb_1} f(s_1) \ \kappa(ds_1|s_2),$ for each bounded continuous function $f : \Sb_1 \to \Rb.$  In view of \eqref{eq:notation:stochastic_kernel}, $\int_{\Omega} f(\phi(s_2,\omega)) \ p(d\omega)=\int_{\Sb_1} f(s_1) \ \kappa(ds_1|s_2),$ and, therefore, statement~(a) holds. For statement (b), stochastic continuity, which is another name for continuity in probability, is a stronger property than continuity in distribution \cite[p. 30]{wentzell_1996_course}. \end{proof}

For the next lemma and the proof of Theorem~\ref{thm:kernels:total_variation}, we shall use the following notation. For a metric space $\Sb,$ we write  $d(x,A) \defeq \inf_{a \in A} \rho_\Sb(x,a)$ for the distance from the point $x \in \Sb$ to the set $A \subset \Sb.$ Given two closed and nonempty subsets $A, B \subset \Sb$ we write the Hausdorff metric
\begin{equation} \label{eq:hausdorff_metric}
    d_H(A, B) \defeq \max\set*{
        \sup_{a \in A} d(a,B),
        \sup_{b \in B} d(b,A)
    },
\end{equation}
which is  possibly infinite. Let $\{C^{(k)}\}_{k=1}^{\infty}$ be a sequence of closed and nonempty subsets of $\Sb.$ We recall that $C^{(k)}$ converges in Hausdorff metric to a closed and nonempty set $C\in\Sb $ if $d_H(C^{(k)}, C) \to 0$ as $k \to \infty,$ written $C^{(k)} \hto C.$


 We say that  a sequence  $\set{C^{(k)}}_{k=1}^{\infty}$ of Borel subsets of $\Rb^n$ converges in the sense of Fr\'echet-Nikodym to a Borel set $C\subset \Rb^n$ if $\lim_{k\to\infty}\lambda^{[n]}(C^{(k)} \triangle C) = 0,$ where $\triangle$ denotes symmetric difference.
 This convergence is equivalent  to convergence in Lebesgue measure of the indicator functions $\oneAlt{C^{(k)}}$ to $\oneAlt{C};$ see~\citet{beer_hausdorff_1974} for an exposition on this topic.
In general, $C^{(k)} \hto C$ neither implies nor is implied by  $\lambda^{[n]}(C^{(k)} \triangle C) \to 0$ even when $C$ and $C^{(k)}$, $k=1,2\ldots,$ are compact sets in $\Rb^n.$ Consider the following examples: if we enumerate $\Qb \cap [0,1] = \set{r_1, r_2, \dots},$ then the sets $C^{(k)} = \set{r_1, r_2, \dots, r_k}$ converge in Hausdorff metric to $[0,1],$ but $\lambda^{[1]}( C^{(k)} \triangle [0,1]) = 1$ for all $k=1,2,\ldots;$ conversely, $\lim_{k\to\infty}\lambda^{[2]}(\tilde{C}^{(k)} \triangle \tilde{C}) = 0$ for $\tilde{C}^{(k)}=[0,k^{-1}]\times [0,1]$ and $\tilde{C}=\{0\}\times [0,0.5],$ but $\lim_{k\to\infty} d_H(\tilde{C}^{(k)},\tilde{C})=0.5.$ The following lemma gives a special case when convergences in Hausdorff and  Fr\'echet-Nikodym senses hold simultaneously; namely, when the sets in question are the images of a function $\phi : \Sb_2 \times \Rb^n \to \Rb^n$ carrying $(s_2,\omega)$ to $\phi(s_2,\omega)$ metrices and satisfying the Diffeomorphic Condition.

  In a metric space $\Sb,$ for  $r>0$ and a set $\mathcal{D}$ we denote $\mathcal{D}^r:=\{x\,:\,d_H(x,D) < r\}.$ If $\mathcal{D}$ is compact and $\mathcal{O}\supset\mathcal{D}$ is its neighborhood,   then there exists $\ell>0$ such that $\mathcal{D}^\ell\subset \mathcal{O}$ \cite[p. 177]{munkrestopology}.
In the proofs of Lemma~\ref{lem:proofs:images_converge_in_measure} and Theorems~\ref{thm:kernels:total_variation} we use  notations $\phi_{s_2}(\omega) \defeq \phi(s_2,\omega)$ and $\mathcal{D}_{s_2} \defeq \phi_{s_2}(\mathcal{D}),$ where    $s_2 \in {\Sb_2},$  $\omega\in \Rb^n,$ and $\mathcal{D}\subset\Rb^n.$  Note that  throughout this paper
$\mathcal{D}_{s_2}^r :=(\mathcal{D}_{s_2})^r.$ 

\begin{lemma} \label{lem:proofs:images_converge_in_measure}
    Let $\Sb_2$ be a metric {space}, and let $\Omega \subset \Rb^n$ be open.
    If a function $\phi : \Sb_2 \times \Rb^n \to \Rb^n$ satisfies the Diffeomorphic Condition, then for every compact set $K \subset \Omega$ and $s_2' \to s_2,$ the convergence of $\phi(s_2', K)$ to  $\phi(s_2, K)$   holds both in Hausdorff and Fr\'echet-Nikodym senses; 
     that is, 
    \begin{equation}\label{eq9.2}
        \lim_{s_2' \to s_2} d_H( \phi(s_2', K), \phi(s_2, K)) =
        \lim_{s_2' \to s_2} \lambda^{[n]}( \phi(s_2', K) \triangle \phi(s_2, K)) = 0.
    \end{equation}
\end{lemma}

\begin{proof}
We note that since $\phi_{s_2}$ is continuous, each $K_{s_2}$ is compact. \citet[Theorem 4.26 and Example~4.13]{rockafellar_variational_2009} imply that if $J : \Rb^{n+1} \to \Rb^{n}$ 
is a continuous mapping and $\set{C^{(k)}}_{k=1}^{\infty} \subset \Rb^{n+1}$ is a uniformly bounded sequence of compact sets such that $C^{(k)} \hto C,$ where $C$ is closed, then $J(C^{(k)}) \hto J(C).$ Now,  let $s_2^{(k)} \to s_2.$ We set
\[
J(\lambda,\omega):= \left\{
\begin{array}{ll}
  \phi(s_2^{(k)},\omega),& \lambda \in \{\frac1k:k=1,2,\ldots\},\,\omega\in K; \\
  \phi(s_2,\omega), & \lambda = 0,\,\omega\in K.
\end{array}
\right.
\]
Note that $J$ is continuous on a compact subset $\{0,\frac1k:k=1,2,\ldots\}\times K$ of $\Rb^{n+1}.$ The Tietze extension theorem allows continuous extension $\tilde{J} : \Rb^{n+1} \to \Rb^{n}$ of $J$ to $\Rb^{n+1}.$ Furthermore $\set{\frac1k} \hto \set{0}$ as $k \to \infty;$ hence $\set{\frac1k}\times K \hto \set{0}\times K.$ Therefore, we may apply~\cite[Theorem 4.26]{rockafellar_variational_2009} for $\tilde{J}$ to conclude that $K_{s_2^{(k)}} := \phi(s_2^{(k)}, K)  = \tilde{J}(\frac1k, K) \hto  \tilde{J}(0,K) = \phi(s_2, K) = K_{s_2}$ as $k \to \infty,$ which establishes convergence in Hausdorff metric. 

 Let us provide a short direct proof of the second equality in \eqref{eq9.2}, which does not use the criterion for $\lim_{k\to\infty}\lambda^{[n]}(C^{(k)} \triangle C) = 0$  from Beer \cite{beer_hausdorff_1974}; see also \cite{lucchetti2020class}.
%
      Let us fix an arbitrary $s_2 \in \Sb_2.$  Before  proving the second equality in \eqref{eq9.2}, we need to establish auxiliary inequality~\eqref{eq:newR1} below. Choose $\ell>0$ such that ${\rm cl}K^\ell\subset \Omega,$ where ${\rm cl}K^\ell$ denotes the closure of $K^\ell.$  This choice is possible in view of the explanations preceding Lemma~\ref{lem:proofs:images_converge_in_measure}. The set ${\rm cl}K^\ell$ is compact because it is a bounded and closed subset of Euclidean space. We apply Berge's maximum theorem to the continuous function $(\hat{s}_2,\omega) \mapsto \|D_\omega \phi_{\hat{s}_2}(\omega)\|^2$ and  compact-valued continuous multifunction $\hat{s}_2\mapsto {\rm cl}K^\ell\ne\emptyset,$ which does not depend on $\hat{s}_2;$ see \cite[Theorem~7, p.~112]{berge1877topological} or \cite[Theorem~1.4]{feinberg_berge_2013}. Therefore, the function $\hat{s}_2\mapsto \max\limits_{\omega\in {\rm cl}K^\ell} \|D_\omega \phi_{\hat{s}_2}(\omega)\|^2$ is continuous. This implies the existence of $M,\delta>0$ such that
    \begin{equation}\label{eq:newR1}
      \|D_\omega \phi_{\hat{s}_2}(\omega)\|\le M<+\infty \mbox{ if } \rho_{\Sb_2}(\hat{s}_2,s_2) \le \delta,\qquad\omega\in {\rm cl}K^\ell.
    \end{equation} 

    Let us continue verifying the second equality in \eqref{eq9.2}. Fix an arbitrary $\varepsilon > 0.$
    Then there exists $r_{s_2}\in(0,\ell)$ such that
    \begin{equation} \label{eq:proofs:hausdorff_eq1}
        \max\set{\lambda^{[n]}(K^r \setminus K), \lambda^{[n]}(K_{s_2}^r \setminus K_{s_2})} < \varepsilon, \qquad r \in [0, r_{s_2}].
    \end{equation}
    because $K^r$ decreases to $K$ and $K_{s_2}^r$ decreases to $K_{s_2}$ as
    $r\downarrow 0.$ 

    Now, let $s_2^{(k)} \to s_2.$ Let us show that there exists $r\in(0,r_{s_2})$ such that
    \begin{equation}\label{eq:IniR}
      K_{s_2^{(k)}}^r\subset \phi_{s^{(k)}_2}(K^{r_{s_2}})\mbox{ if }\rho_{\Sb_2}(s_2^{(k)},s_2) \leq \delta,
    \end{equation}
    where $\delta>0$ is defined in \eqref{eq:newR1}. Indeed, let $S_2 := \{s_2, s_2^{(k)} : k = k_1, k_1 + 1, \ldots\},$ where $k_1$ is the smallest
    natural
    number such that $S_2$ is contained within the open ball $B^\delta(s_2)$ of radius $\delta$ centered at $s_2.$   Consider the compact metric spaces $(S_2,\rho_{\Sb_2})$ and $\tilde{\Sb}:=S_2\times{\rm cl}K^\ell.$ The mapping $(\hat{s}_2,\omega)\mapsto \tilde{\phi}(\hat{s}_2,\omega):=(\hat{s}_2,\phi_{\hat{s}_2}(\omega))$ acting from $\tilde{\Sb}$ onto the metric space $(\tilde{\phi}(\tilde{\Sb}),\rho_{\Sb_2\times\Rb^n})$ is bijective and continuous, where continuity follows from continuity of composition of continuous functions. 
    In view of \cite[Theorem~26.6]{munkrestopology}, since $\tilde{\Sb}$ is a compact metric space, the mapping $\tilde{\phi}:\tilde{\Sb}\to \Phi(\tilde{\Sb})$ is a homeomorphism. Therefore, the image $\tilde{\phi}(S_2 \times K^{r_{s_2}})$ of the open set $S_2 \times K^{r_{s_2}}$ is open in $\tilde{\phi}(\tilde{\Sb}).$ In addition,  $\tilde{\phi}(S_2 \times K)\subset \tilde{\phi}(S_2 \times K^{r_{s_2}}),$  where $\tilde{\phi}(S_2 \times K)$ is a compact set in $\tilde{\phi}(\tilde{\Sb})$  as a continuous image of a compact set.  As explained in the paragraph preceding Lemma~\ref{lem:proofs:images_converge_in_measure}, there exists $r\in (0,r_{s_2})$ such that
    $
    (\tilde{\phi}(S_2 \times K))^r \subset \tilde{\phi}(S_2 \times K^{r_{s_2}}),
    $
    which for $k\ge k_1$ implies
    $
    \{s^{(k)}_2\}\times K^r_{s^{(k)}_2}\subset (\{s^{(k)}_2\}\times K_{s^{(k)}})^r\subset (\tilde{\phi}(S_2\times K))^r\subset\tilde{\phi}(S_2\times K^{r_{s_2}}),
    $
    which in its turn implies \eqref{eq:IniR}.

    Since $K_{s_2^{(k)}} \hto K_{s_2},$ there exists $k_2= k_1,k_1+1,\ldots$ such that $K_{s_2^{(k)}} \subset K_{s_2}^{r}$  and $K_{s_2} \subset K_{s_2^{(k)}}^{r}$ for each $k=k_2,k_2+1,\ldots.$ Hence,  for each  $k=k_2,k_2+1,\ldots,$ by \eqref{eq:proofs:hausdorff_eq1}, $\lambda^{[n]}(K_{s_2^{(k)}} \setminus K_{s_2})   < \varepsilon,$ 
    and, by \eqref{eq:newR1}--\eqref{eq:IniR},
    \[
    \lambda^{[n]}(K_{s_2} \setminus K_{s_2^{(k)}}) \leq\lambda^{[n]}( \phi_{s^{(k)}_2}(K^{r_{s_2}}) \setminus K_{s_2^{(k)}}) = \int_{K^{r_{s_2}} \setminus K} \abs{\det D_\omega\phi_{s_2^{(k)}}(\omega)} \ d\omega \leq M \cdot \lambda^{[n]}(K^{r_{s_2}} \setminus K)<M\varepsilon,
    \]
     where the first inequality holds because $K_{s_2} \subset K_{s_2^{(k)}}^{r}\subset\phi_{s^{(k)}_2}(K^{r_{s_2}}),$ the equality holds following the change-of-variables formula; see, e.g.~\cite[Theorem 7.26]{walter_rudin_1966}, the second inequality follows from \eqref{eq:newR1}, and the last inequality follows from \eqref{eq:proofs:hausdorff_eq1}. Therefore, $\mathop{\rm lim\,sup}\limits_{k\to \infty} \lambda^{[n]}(K_{s_2^{(k)}} \triangle K_{s_2}) < (1+M)\varepsilon.$
  Since $\varepsilon > 0$ is arbitrary, $\lim\limits_{k\to\infty} \lambda^{[n]}( K_{s^{(k)}_2} \triangle K_{s_2}) = 0.$
\end{proof}



\begin{proof}[{Proof of Theorem~\ref{thm:kernels:total_variation}}]
    We first observe that continuity of $\phi$  in total variation {\it wrt} $p$  is necessary and sufficient for continuity in total variation of $\kappa.$ Indeed, since for all $s_2, s_2' \in \Sb_2$
    \begin{equation*}
        \sup_{B \in \Bc(\Sb_1)} \abs{
            \kappa(B|s_2') - \kappa(B|s_2)
        } =
        \sup_{B \in \Bc(\Sb_1)} \abs*{
            \int_{\Omega} \mathbf{1}\{\phi(s_2',\omega) \in B\} -
            \mathbf{1}\{\phi(s_2, \omega) \in B\} \
            p(d\omega)
        },
    \end{equation*}
    $\phi$ is continuous in total variation {\it wrt} $p$ if and only if $\kappa$ is continuous in total variation  that is, statement~(a) holds.

    Now,  let us  prove statement~(b). Since $p \ll \lambda^{[n]},$ let $f:\Omega\to\Rb$ be a nonnegative Borel measurable function such that $dp = f \ d\lambda^{[n]}.$ Let $s_2^{(k)} \to s_2.$ We need to show that
    \begin{equation}
        \adjustlimits\lim_{s_2^{(k)} \to s_2} \sup_{B \in \Bc(\Sb_2)} \abs{\kappa(B|s_2^{(k)}) - \kappa(B|s_2)} = 0.
    \end{equation}
    Let $\varepsilon > 0$ be arbitrary, and let $B \in \Bc(\Omega).$ Since $\int_{\Omega} f(\omega) \ \lambda^{[n]}(d\omega) =p(\Omega)= 1,$
    by Markov's inequality $p(A) \geq 1-\frac{\varepsilon}{3}$ for   the set $A: = \set{\omega \in \Omega : 0 \le f(\omega) \leq \frac{3}{\varepsilon}}.$ 
     Since $p$ is a regular probability measure, 
    we can fix a compact set $K \subset A$ such that $p(K) \geq p(A) - \frac{\varepsilon}{3}.$ On this set $K,$ it follows that 
    $f \leq M$ for some $M > 0,$ and furthermore $p(K) \geq 1-\frac{2\varepsilon}{3}.$
    Lusin's theorem implies the existence of a compact subset $C\subset K$ with $p(K\setminus C) < \frac{\varepsilon}3$ such that the restriction of
$f$  to $C$ is continuous. Observe that $p(C) \geq 1-\varepsilon.$

We now show that the functions $\abs{\det D_\omega \phi_{s_2}}$ and $\abs{\det D_\omega\phi_{s_2^{(k)}}}$ are uniformly bounded away from $0$ on $C.$ The set $\tilde{K}: = (\set{s_2} \cup \set{s_2^{(k)}}_{k=1}^\infty) \times C$ is compact. By  assumptions,  the function $(s_2,\omega) \mapsto \abs{\det D_\omega \phi_{s_2}(\omega)}$ is continuous and strictly positive on the compact set $\tilde{K},$  and we can fix $\beta > 0$ such that $\abs{\det D_\omega \phi_{s_2}(\omega)} \geq \beta$ for all $(s_2,\omega) \in \tilde{K}.$ Now, let us use the representation $\Omega=C\cup(\Omega\setminus C)$ and then use change-of-coordinates to obtain
    \begin{align*}
       & \abs*{\kappa(B|s_2^{(k)}) - \kappa(B|s_2)}  =
        \abs*{
            \int_{\Omega} \mathbf{1}\{\phi_{s_2^{(k)}}(\omega) \in B\} \ p(d\omega) -
            \int_{\Omega} \mathbf{1}\{\phi_{s_2}(\omega) \in B\} \ p(d\omega)}
        \\ & \leq \varepsilon +
        \abs*{
            \int_{C} \mathbf{1}\{\phi_{s_2^{(k)}}(\omega) \in B\} \ f(\omega)\, d\omega -
            \int_{C} \mathbf{1}\{\phi_{s_2}(\omega) \in B\} \ f(\omega)\, d\omega
        } \\
        &  =  \varepsilon
        +\abs*{
            \int_{B \cap \phi_{s_2^{(k)}}}(C) f(\phi^{-1}_{s_2^{(k)}}(s_1))
            \abs{\det D_{s_1} \phi^{-1}_{s_2^{(k)}}(s_1)}
            \ ds_1 -
            \int_{B \cap \phi_{s_2}(C)} f(\phi^{-1}_{s_2}(s_1))
            \abs{\det {D_{s_1}}\phi^{-1}_{s_2}(s_1)}
            \ ds_1
        }.
    \end{align*}
    The inequalities $\abs{\det D_{\omega} \phi_{\hat{s}_2}(\omega)} \geq \beta$ on $\tilde{K}$  and $0\le f \le M$ on $C$ imply that \[0\le f(\phi^{-1}_{\hat{s}_2}(s_1))\abs{\det D_{s_1} \phi^{-1}_{\hat{s}_2}(s_1)} \leq M\beta\inv\] on the  set $  \bigcup_{k=0}^{\infty} \{s_2^{(k)}\}\times\phi_{s_2^{(k)}}(C),$ where $s_2^{(0)}:=s_2.$ 
    Therefore, we obtain the estimate
    \begin{align}
        \begin{split}
            \abs*{\kappa(B|s_2^{(k)}) - \kappa(B|s_2)} &\leq
            \varepsilon + \beta\inv M\cdot \lambda^{[n]}(B \cap [\phi_{s_2^{(k)}} (C) \triangle \phi_{s_2} (C)]) \\
             +
            \int_{B \cap \phi_{s_2^{(k)}} (C) \cap \phi_{s_2} (C)} & \abs*{f(\phi^{-1}_{s_2^{(k)}}(s_1))\abs{\det D_{s_1} \phi^{-1}_{s_2^{(k)}}(s_1)} - f(\phi^{-1}_{s_2}(s_1))
             \abs{\det D_{s_1} \phi^{-1}_{s_2}(s_1)}} \ ds_1 \\
            &\leq\label{eq:Shaffe}
            \varepsilon + \beta\inv M\cdot \lambda^{[n]}([\phi_{s_2^{(k)}} (C) \triangle \phi_{s_2} (C)]) \\
             +
            \int_{\phi_{s_2^{(k)}} (C) \cap \phi_{s_2} (C)} & \abs*{f(\phi^{-1}_{s_2^{(k)}}(s_1))\abs{\det D_{s_1} \phi^{-1}_{s_2^{(k)}}(s_1)} - f(\phi^{-1}_{s_2}(s_1))
             \abs{\det D_{s_1} \phi^{-1}_{s_2}(s_1)}} \ ds_1,
        \end{split}
    \end{align}
    and we observe that this estimate does not depend on the choice of $B.$

   Lemma~\ref{lem:proofs:images_converge_in_measure} implies that $\lambda^{[n]}(\phi_{s_2^{(k)}}(C) \triangle \phi_{s_2}( C)) \to 0$ as $k\to\infty.$ Let us prove that
    \begin{equation} \label{eq:proofs:intersection}
        \lim_{s_2^{(k)} \to s_2} \int_{\phi_{s_2^{(k)}} (C) \cap \phi_{s_2} (C)} \abs*{f(\phi^{-1}_{s_2^{(k)}}(s_1))\abs{\det D_{s_1} \phi^{-1}_{s_2^{(k)}}(s_1)} - f(\phi^{-1}_{s_2}(s_1))
             \abs{\det D_{s_1} \phi^{-1}_{s_2}(s_1)}} \ ds_1 = 0.
    \end{equation}
   For this purpose we firstly note that
    \begin{equation}\label{eq:new1}
        \abs*{f(\phi^{-1}_{s_2^{(k)}}(s_1))\abs{\det {D_{s_1}} \phi^{-1}_{s_2^{(k)}}(s_1)} - f(\phi^{-1}_{s_2}(s_1))
             \abs{\det {D_{s_1}} \phi^{-1}_{s_2}(s_1)}}\leq M\beta\inv,
    \end{equation}
{for each $k.$ Secondly, we note that
    the pointwise convergence
        \begin{equation}\label{eq:new2}
          \begin{aligned}
      I^{(k)}(s_1) & :=\mathbf{1}\{s_1 \in \phi_{s_2^{(k)}} (C) \cap \phi_{s_2} (C)\} \times\\
      \times & \abs*{f(\phi^{-1}_{s_2^{(k)}}(s_1))\abs{\det {D_{s_1}} \phi^{-1}_{s_2^{(k)}}(s_1)} - f(\phi^{-1}_{s_2}(s_1))
             \abs{\det {D_{s_1}} \phi^{-1}_{s_2}(s_1)}}\to 0
    \end{aligned}
        \end{equation}
                  holds. Indeed, $I^{(k)}(s_1) \le I_1^{(k)}(s_1) + I_2^{(k)}(s_1),$ where
    \begin{align*}
        I_1^{(k)}(s_1) & := \mathbf{1}\{s_1 \in \phi_{s_2^{(k)}} (C) \cap \phi_{s_2} (C)\}\abs{f(\phi^{-1}_{s_2^{(k)}}(s_1))}\,\abs*{\abs{\det {D_{s_1}} \phi^{-1}_{s_2^{(k)}}(s_1)} -
             \abs{\det {D_{s_1}} \phi^{-1}_{s_2}(s_1)}}\\ &\le M {\mathbf{1}\{s_1 \in \phi_{s_2^{(k)}} (C) \cap \phi_{s_2} (C)\}}\abs*{\abs{\det {D_{s_1}} \phi^{-1}_{s_2^{(k)}}(s_1)} - \abs{\det {D_{s_1}} \phi^{-1}_{s_2}(s_1)}}\to 0 \quad {\rm as}\ k\to\infty,
    \end{align*}
     where the inequality holds because $0\le f\le M$ on $C,$ and the convergence follows from Remark~\ref{rem:new}; and
     \begin{align*}
        I_2^{(k)}(s_1) & := \mathbf{1}\{s_1 \in \phi_{s_2^{(k)}} (C) \cap \phi_{s_2} (C)\}\abs{\det {D_{s_1}} \phi^{-1}_{s_2}(s_1)}\cdot \abs{f(\phi^{-1}_{s_2^{(k)}}(s_1)) - f(\phi^{-1}_{s_2}(s_1))
             }\\
             & =\mathbf{1}\{\phi^{-1}_{s_2^{(k)}}(s_1) \in C \}\mathbf{1}\{\phi^{-1}_{s_2}(s_1) \in C \} \abs{\det {D_{s_1}} \phi^{-1}_{s_2}(s_1)} \cdot \abs{f(\phi^{-1}_{s_2^{(k)}}(s_1)) - f(\phi^{-1}_{s_2}(s_1))
             }\to 0\quad {\rm as}\ k\to\infty
    \end{align*}
                because, according to Remark~\ref{rem:new}, $\phi^{-1}_{s_2^{(k)}}(s_1) \to \phi^{-1}_{s_2}(s_1),$  when $\phi^{-1}_{s_2^{(k)}}(s_1) \in C\subset\Omega,$ and $f$ restricted to $C$ is continuous.
             Then \eqref{eq:new1}, \eqref{eq:new2}, and} 
             the dominated convergence theorem implies~\eqref{eq:proofs:intersection}. Therefore,
    \begin{equation*}
        \mathop{\rm lim\, sup}\limits_{s_2^{(k)} \to s_2} \sup_{B \in \Bc(\Sb_2)} \abs*{\kappa(B|s_2^{(k)}) - \kappa(B|s_2)} \leq \varepsilon,
    \end{equation*}
    and since $\varepsilon > 0$ was arbitrary, it follows that $\kappa$ is continuous in total variation, that is, statement~(b) holds.
\end{proof}


\begin{remark}\label{rem:Shaffe}
The proof of Theorem~\ref{thm:kernels:total_variation}(b) can be  simplified under additional assumptions.  The examples of such assumptions are the case $n=1$ considered in
\cite[Section 8.1]{feinberg_partially_2016} and the case when $p$ has a continuous density and the  images $\phi(s_2,\Omega)$ are the same for all $s_2\in\Sb_2$ considered in  in  \cite[Example 1.3.2(v)]{KY}. In the general case,  convergence to zero of the  term $\lambda^{[n]}(B \cap [\phi_{s_2^{(k)}} (C) \triangle \phi_{s_2} (C)])$ in \eqref{eq:Shaffe}, which is proved by using Lemma~\ref{lem:proofs:images_converge_in_measure}, plays the key role in the proof of Theorem~\ref{thm:kernels:total_variation}(b).
In particular,  \cite[Example 1.3.2(v)]{KY} may create an impression that, if the probability $p$ has a continuous density, then Theorem~\ref{thm:kernels:total_variation}(b) can be proved easily by applying Scheff\'e's theorem.  However, the proof in  \cite[Example 1.3.2(v)]{KY} uses a hidden assumption that the images $\phi(s_2,\Omega)$ do not depend on $s_2.$ \end{remark}

The next lemma is needed to prove Theorem~\ref{thm:feller:filter_kernel_continuous}. It describes the relationships between the statements of Assumption~\ref{ass:filters:assumption_1} and~\ref{ass:filters:assumption_2} to continuity properties of the stochastic kernels $\Tc$ and $Q$ defined in \eqref{eq:filters:model_kernels}.

\begin{lemma} \label{lem:proofs:kernel_assumptions}
    Consider the stochastic kernels $\Tc$ on $\Xb$ given $\Xb \times \Ab$ and $Q$ on $\Yb$ given $\Ab \times \Xb$ defined in~\eqref{eq:filters:model_kernels}. The following statements hold:
    \begin{enumerate}
        \item[\rm(a)] Under Assumption~\ref{ass:filters:assumption_1}(i), the stochastic kernel $\Tc$ is weakly continuous, and under Assumption~\ref{ass:filters:assumption_2}(i) the stochastic kernel $\Tc$ is continuous in total variation. 
        \item[\rm(b)] Under Assumption~\ref{ass:filters:assumption_1}(ii), the stochastic kernel $Q$ is continuous in total variation, and under Assumption~\ref{ass:filters:assumption_2}(ii), the stochastic kernel $Q$ is continuous in $a$ in total variation. 
    \end{enumerate}
\end{lemma}

\begin{proof}
    The lemma follows directly from Theorems~\ref{thm:kernels:weak} and \ref{thm:kernels:total_variation}. In particular, Assumption~\ref{ass:filters:assumption_1}(i) and Theorem~\ref{thm:kernels:weak} imply that $\Tc$ is weakly continuous, whereas Assumption~\ref{ass:filters:assumption_2}(i) and Theorem~\ref{thm:kernels:total_variation} imply that $\Tc$ is continuous in total variation. Similarly, Assumption~\ref{ass:filters:assumption_1}(ii) and Theorem~\ref{thm:kernels:total_variation} imply that $Q$ is continuous in total variation. If $G$ does not depend on $a,$ then continuity in $a$ in total variation is immediate. On the other hand, Assumption~\ref{ass:filters:assumption_2}(ii) and Theorem~\ref{thm:kernels:total_variation} imply that $Q$ is continuous in $a$ in total variation if $G$ depends on $a.$
\end{proof}


\begin{proof}[{Proof of Theorem~\ref{thm:feller:filter_kernel_continuous}}]
    Lemma~\ref{lem:proofs:kernel_assumptions} and Theorem~\ref{thm:feller:semi-uniform_feller} imply that Assumptions~\ref{ass:filters:assumption_1} and~\ref{ass:filters:assumption_2} are each sufficient for $P$ to be semi-uniform Feller. 
\end{proof}


\begin{proof}[Proof of Theorem~\ref{thm:filters:main_result}]
    The conclusions follow from Theorems~\ref{thm:feller:filter_kernel_continuous} and~\ref{thm:feller:semi-uniform_Feller_weak_continuous_filter}.
\end{proof}


\begin{proof}[{Proof of Theorem~\ref{thm:pomdp1s:results}}]
    Under the given conditions, Assumption~\ref{ass:filters:assumption_1} is satisfied for the \pomdpOne{}. Hence, Lemma~\ref{lem:proofs:kernel_assumptions} implies that $\Tc$ is weakly continuous and $Q_1$ is continuous in total variation, and Theorem~\ref{thm:pomdp1s:semi-uniform_feller} implies that $P$ is semi-uniform Feller.
\end{proof}

Finally, we close this section by describing an alternative proof of Corollary~\ref{cor:kernels:aumann_multidim} for the case $\Sb_1 = \Rb^n.$ The construction in the following proof enables a straightforward verification of the Diffeomorphic Condition.

\begin{proof}[Alternative proof of Corollary~\ref{cor:kernels:aumann_multidim} when $\Sb_1 = \Rb^n$]
    Let $A_1, \dots, A_n \in \Bc(\Rb)$ be arbitrary.
        According to~\citet[Proposition 7.27]{bertsekas_stochastic_1996}, the stochastic kernel $\kappa$ can be decomposed on products of the form $A_1 \times \dots \times A_n$ by the conditional representation
        \begin{equation} \label{eq:proofs:conditional_decomposition}
            \kappa(A_1 \times \cdots \times A_n|s_2) =
            \int_{A_1} \cdots \int_{A_n} \kappa_{n-1}(dx_n|s_2,x_1,\dots,x_{n-1}) \cdots \kappa_1(dx_2|s_2,x_1) \ \kappa_0(dx_1|s_2),
        \end{equation}
        where $\kappa_0$ is a stochastic kernel on $\Rb$ given $\Sb_2,$ and for each $j=1, 2, \dots, n-1,$ $\kappa_j$ is a stochastic kernel on $\Rb$ given $\Sb_2 \times \Rb^{j}.$ In view of Theorem~\ref{thm:kernels:aumann}, there exist Borel measurable functions $\phi_0 : \Sb_2 \times [0,1] \to \Rb$ and $\phi_j : (\Sb_2 \times \Rb^j) \times [0,1] \to \Rb$ ($j=1,\dots,n-1$) such that for each $B \in \Bc(\Rb)$ and $s_2 \in \Sb_2,$ there holds
        \begin{align*}
            \kappa_0(B|s_2) &= \int_0^1 \mathbf{1}\{\phi_0(s_2,\omega_1) \in B\} \ d\omega_1 \\
            \kappa_j(B|s_2,x_1,\dots,x_j) &=
            \int_0^1 \mathbf{1}\{\phi_j(s_2,x_1,\dots,x_j,\omega_{j+1}) \in B\} \ d\omega_{j+1}, && x_1, \dots, x_j \in \Rb.
        \end{align*}
        We now consider the following construction. Given $s_2 \in \Sb_2$ and $(\omega_1, \dots, \omega_n) \in [0,1]^n,$ let $x_1 = \phi_0(s_2,\omega_1).$ Proceeding by induction, if we have obtained $x_1, \dots, x_{j-1},$ we calculate $x_j = \phi_{j-1}(s_2, x_1, \dots, x_{j-1}, \omega_j).$ Define the function $\phi : \Sb_2 \times [0,1]^n \to \Rb^n$ by the formula $\phi(s_2, \omega) = (x_1, \dots, x_n),$ where each $x_j$ ($j=1,\dots,n$) is calculated from the above procedure. As a composition of Borel-measurable functions, $\phi$ is a Borel measurable function. Furthermore, from~\eqref{eq:proofs:conditional_decomposition} it follows that
        \begin{equation*}
            \kappa(A_1 \times \cdots \times A_n|s_2) = \int_{[0,1]^n} \mathbf{1}\{\phi(s_2,\omega) \in A_1 \times \dots \times A_n\} \ d\omega,
        \end{equation*}
        and therefore \eqref{eq:kernels:multidim} holds. 
\end{proof}

    The benefit of the above proof is that it can be used to construct functions $\phi$ for which the representation~\eqref{eq:kernels:multidim} holds, and the Diffeomorphic Condition is satisfied. Let us consider a simple example involving a multivariate normal distribution. Let $\Sb_2 \defeq (-1,1),$ and consider the stochastic kernel $\kappa$ on $\Rb^2$ given $\Sb_2$ such that $\kappa(\:\cdot\:|s_2)$ is centered 
    at the origin with covariance matrix $\Sigma = \begin{bsmallmatrix}
    1 & s_2\\
    s_2 & 1
\end{bsmallmatrix}$ for $s_2 \in \Sb_2;$ that is,
\begin{align*}
\kappa(B|s_2) &\defeq \frac{1}{2\pi \sqrt{(1-s_2^2)}}
\int_B
\exp(-\tfrac{1}{2}x^T\Sigma^{-1} x) \ dx,
&&
B \in \Bc(\Rb^2), \quad s_2 \in \Sb_2,
\end{align*}
where $x = (x_1, x_2)^T.$ 
The parameter $s_2$ represents the linear correlation between the components of $x.$ We note that the function $\phi$ provided in the first proof of Corollary~\ref{cor:kernels:aumann_multidim} at the beginning of this section does not satisfy the Diffeomorphic Condition, since the components of $\phi$ in that construction are identical, which implies that $D_\omega \phi$ is singular. However, by following the construction in the alternative proof, we shall recover a function $\phi$ that satisfies the Diffeomorphic Condition.
To construct $\phi$ for the representation~\eqref{eq:kernels:multidim}, we observe that the kernel $\kappa$ has the conditional representation
\begin{align*}
    \kappa(A_1 \times A_2 | s_2) &=
    \int_{A_1} \int_{A_2} \kappa_1(dx_2|s_2, x_1) \ \kappa_0(dx_1|s_2),
    &&
    A_1, A_2 \in \Bc(\Rb),
\end{align*}
where $\kappa_0$ is a stochastic kernel on $\Rb$ given $\Sb_2$ and $\kappa_1$ is a stochastic kernel on $\Rb$ given $\Sb_2 \times \Rb.$ In this instance, the marginal and conditional distributions for the multivariate normal distribution are well-known, so we may use
\begin{align*}
    \kappa_0(A_1 | s_2) &\defeq \frac{1}{\sqrt{2\pi}} \int_{A_1}
    \exp(-\tfrac{x_1^2}{2})\ dx_1, &&
    A_1 \in \Bc(\Rb),  \\
    \kappa_1(A_2 | s_2, x_1) &\defeq  \frac{1}{\sqrt{2\pi(1-s_2^2)}} \int_{A_2}
    \exp(-\tfrac{(x_2-s_2 x_1)^2}{2(1-s_2^2)}) \ dx_2, &&
    A_2 \in \Bc(\Rb).
\end{align*}
In this form, both $\kappa_0$ and $\kappa_1$ are distributions of univariate 
normal random variables. According to Theorem~\ref{thm:kernels:aumann}, we can fix functions $\phi_0 : \Sb_2 \times [0,1] \to \Rb$ and $\phi_1 : \Sb_2 \times \Rb \times [0,1] \to \Rb$ such that~\eqref{eq:kernels:kernel_representation} holds for $\kappa_0$ and $\kappa_1,$ respectively. In particular, letting $\Phi$ denote the cumulative distribution function for a standard normal distribution, it follows that  for $(x_1,x_2)\in\Rb^2$
\begin{align*}
    \phi_0(s_2, \omega_1) &\defeq
    \Phi^{-1}(\omega_1) \\
    \phi_1(s_2, x_1, \omega_2) &\defeq
    s_2x_1 + \sqrt{(1-s_2^2)} \Phi^{-1}(\omega_2)
\end{align*}
are such representative functions. Finally, for $\omega = (\omega_1, \omega_2),$ letting
\begin{equation*}
    \phi(s_2, \omega) \defeq
    \begin{bmatrix}
        \Phi^{-1}(\omega_1) \\
       s_2\Phi^{-1}(\omega_1) + \sqrt{(1-s_2^2)} \Phi^{-1}(\omega_2)
    \end{bmatrix},
\end{equation*}
it follows that the representation~\eqref{eq:kernels:multidim} holds, and $\phi$ is continuous. To confirm that $\phi$ satisfies the Diffeomorphic Condition, we evaluate the determinant of $D_\omega \phi:$
\begin{equation*}
    \det D_\omega \phi =
    \begin{vmatrix}
        [\Phi'(\Phi^{-1}(\omega_1))]^{-1} & 0 \\
        s_2[\Phi'(\Phi^{-1}(\omega_1))]^{-1} & \sqrt{(1-s_2^2)}[\Phi'(\Phi^{-1}(\omega_2))]^{-1}
    \end{vmatrix} =
    \frac{\sqrt{1-s_2^2}}{\Phi'(\Phi^{-1}(\omega_1))\Phi'(\Phi^{-1}(\omega_2))} \ne 0
    \end{equation*}
for all $s_2 \in \Sb_2.$ Hence, $D_\omega \phi$ is nonsingular, so $\phi$ satisfies the Diffeomorphic Condition.

\bibliography{POMDPs_from_equations_FIKK_20240821}

\begin{thebibliography}{43}
\providecommand{\natexlab}[1]{#1}
\providecommand{\url}[1]{\texttt{#1}}
\expandafter\ifx\csname urlstyle\endcsname\relax
  \providecommand{\doi}[1]{doi: #1}\else
  \providecommand{\doi}{doi: \begingroup \urlstyle{rm}\Url}\fi

\bibitem[Aoki(1965)]{aoki_1965}
M.~Aoki.
\newblock Optimal control of partially observable {Markovian} systems.
\newblock \emph{Journal of the Franklin Institute}, 280\penalty0 (5):\penalty0
  367--386, 1965.

\bibitem[{\r{A}}str\"{o}m(1965)]{astrom_1965}
K.~J. {\r{A}}str\"{o}m.
\newblock Optimal control of {Markov} processes with incomplete state
  information.
\newblock \emph{Journal of Mathematical Analysis and Applications}, 10\penalty0
  (1):\penalty0 174--205, 1965.

\bibitem[Aumann(1964)]{aumann_1964_mixed}
R.~J. Aumann.
\newblock Mixed and behavior strategies in infinite exstensive games.
\newblock \emph{Advances in Game Theory}, 52:\penalty0 627--650, 1964.

\bibitem[Beer(1974)]{beer_hausdorff_1974}
G.~A. Beer.
\newblock The {Hausdorff} metric and convergence in measure.
\newblock \emph{Michigan Mathematical Journal}, 21\penalty0 (1):\penalty0
  63--64, 1974.

\bibitem[Bensoussan(2011)]{bensoussan2011dynamic}
A.~Bensoussan.
\newblock \emph{Dynamic Programming and Inventory Control}.
\newblock IOS Press, Amsterdam, 2011.

\bibitem[Berge(1963)]{berge1877topological}
C.~Berge.
\newblock \emph{Topological spaces: Including a treatment of multi-valued
  functions, vector spaces and convexity}.
\newblock Oliver \& Boyd, Edinburgh-London, 1963.

\bibitem[Bertsekas and Shreve(1996)]{bertsekas_stochastic_1996}
D.~Bertsekas and S.~E. Shreve.
\newblock \emph{Stochastic {Optimal} {Control}: the {Discrete}-{Time} {Case}}.
\newblock Athena Scientific, Belmont, 1996.

\bibitem[Downing et~al.(1980)Downing, Pike, and
  Morrison]{downing_1980_application}
D.~J. Downing, D.~H. Pike, and G.~W. Morrison.
\newblock Application of the {Kalman} filter to inventory control.
\newblock \emph{Technometrics}, 22\penalty0 (1):\penalty0 17--22, 1980.

\bibitem[Dynkin(1965)]{dynkin_1965}
E.~B. Dynkin.
\newblock Controlled random sequences.
\newblock \emph{Theory of Probability and Its Applications}, 10\penalty0
  (1):\penalty0 1--14, 1965.

\bibitem[Dynkin and Yushkevich(1979)]{dynkin_1979_controlled}
E.~B. Dynkin and A.~A. Yushkevich.
\newblock \emph{Controlled {Markov} Processes}.
\newblock Springer-Verlag, New York, 1979.

\bibitem[Feinberg(2016)]{feinberg_2016_optimality}
E.~A. Feinberg.
\newblock Optimality conditions for inventory control.
\newblock In A.~Gupta and A.~Capponi, editors, \emph{Tutorials in Operations
  Research: Optimization Challenges in Complex, Networked and Risky Systems},
  pages 14--45. INFORMS, Cantonsville, MD, 2016.

\bibitem[Feinberg and Kasyanov(2021)]{feinberg_2021_mdps}
E.~A. Feinberg and P.~O. Kasyanov.
\newblock {MDP}s with setwise continuous transition probabilities.
\newblock \emph{Operations Research Letters}, 49\penalty0 (5):\penalty0
  734--740, 2021.

\bibitem[Feinberg and Kasyanov(2023)]{feinberg_2023_equivalent}
E.~A. Feinberg and P.~O. Kasyanov.
\newblock Equivalent conditions for weak continuity of nonlinear filters.
\newblock \emph{Systems \& Control Letters}, 173, 2023.

\bibitem[Feinberg and Kraemer(2023)]{feinberg2022continuity}
E.~A. Feinberg and D.~N. Kraemer.
\newblock Continuity of discounted values and the structure of optimal policies
  for periodic-review inventory systems with setup costs.
\newblock \emph{Naval Research Logistics (NRL)}, 70\penalty0 (5):\penalty0
  480--492, 2023.

\bibitem[Feinberg and Lewis(2018)]{FL2018}
E.~A. Feinberg and M.~E. Lewis.
\newblock On the convergence of optimal actions for {M}arkov decision processes
  and the optimality of $(s, {S})$ inventory policies.
\newblock \emph{Naval Research Logistics (NRL)}, 65\penalty0 (8):\penalty0
  619--637, 2018.

\bibitem[Feinberg et~al.(2012)Feinberg, Kasyanov, and
  Zadoianchuk]{feinberg_average_cost_2012}
E.~A. Feinberg, P.~O. Kasyanov, and N.~V. Zadoianchuk.
\newblock Average cost {Markov} decision processes with weakly continuous
  transition probabilities.
\newblock \emph{Mathematics of Operations Research}, 37\penalty0 (4):\penalty0
  591--607, 2012.

\bibitem[Feinberg et~al.(2014)Feinberg, Kasyanov, and
  Voorneveld]{feinberg_berge_2013}
E.~A. Feinberg, P.~O. Kasyanov, and M.~Voorneveld.
\newblock {Berge's} maximum theorem for noncompact image sets.
\newblock \emph{Journal of Mathematical Analysis and Applications},
  413\penalty0 (2):\penalty0 1040--1046, 2014.

\bibitem[Feinberg et~al.(2016)Feinberg, Kasyanov, and
  Zgurovsky]{feinberg_partially_2016}
E.~A. Feinberg, P.~O. Kasyanov, and M.~Z. Zgurovsky.
\newblock Partially observable total-cost {Markov} decision processes with
  weakly continuous transition probabilities.
\newblock \emph{Mathematics of Operations Research}, 41\penalty0 (2):\penalty0
  656--681, 2016.

\bibitem[Feinberg et~al.(2022)Feinberg, Kasyanov, and
  Zgurovsky]{feinberg_markov_2022}
E.~A. Feinberg, P.~O. Kasyanov, and M.~Z. Zgurovsky.
\newblock {Markov} decision processes with incomplete information and
  semi-uniform {F}eller transition probabilities.
\newblock \emph{SIAM Journal on Control and Optimization}, 60\penalty0
  (4):\penalty0 2488--2513, 2022.

\bibitem[Feinberg et~al.(2023)Feinberg, Kasyanov, and
  Zgurovsky]{feinberg_2023_semi}
E.~A. Feinberg, P.~O. Kasyanov, and M.~Z. Zgurovsky.
\newblock Semi-uniform {F}eller stochastic kernels.
\newblock \emph{Journal of Theoretical Probability}, 36:\penalty0 2262--2283,
  2023.

\bibitem[Gihman and Skorohod(1979)]{gihman_1979_controlled}
I.~I. Gihman and A.~V. Skorohod.
\newblock \emph{Controlled {Stochastic} {Processes}}.
\newblock Springer, New York, 1979.

\bibitem[Guo and Mishra(2017)]{guo_constrained_2017}
Y.~Guo and S.~Mishra.
\newblock Constrained optimal iterative learning control with mixed-norm cost
  functions.
\newblock \emph{Mechatronics}, 43:\penalty0 56--65, 2017.

\bibitem[Hendricks et~al.(2008)Hendricks, Jannerup, and
  S\o{}renson]{hendricks_2008_linear}
E.~Hendricks, O.~Jannerup, and P.~H. S\o{}renson.
\newblock \emph{Linear {Systems} {Control}}.
\newblock Springer-Verlag, Berlin, 2008.

\bibitem[Hern\'andez-Lerma(1989)]{hernandez-lerma_adaptive_1989}
O.~Hern\'andez-Lerma.
\newblock \emph{Adaptive {Markov} Control Processes}.
\newblock Springer-Verlag, New York, 1989.

\bibitem[Kara and Yüksel(2020)]{KY}
A.~D. Kara and S.~Yüksel.
\newblock Robustness to incorrect system models in stochastic control.
\newblock \emph{SIAM Journal on Control and Optimization}, 58\penalty0
  (2):\penalty0 1144--1182, 2020.

\bibitem[Kara et~al.(2019)Kara, Saldi, and Y\"uksel]{kara_2019_weak}
A.~D. Kara, N.~Saldi, and S.~Y\"uksel.
\newblock Weak {F}eller property of non-linear filters.
\newblock \emph{Systems and Control Letters}, 134:\penalty0 104512, 2019.

\bibitem[Kifer(1986)]{kifer_1986_ergodic}
Y.~Kifer.
\newblock \emph{Ergodic Theory of Random Transformations}.
\newblock Birkh{\"a}user, Boston, 1986.

\bibitem[Kwakernaak et~al.(1974)Kwakernaak, Sivan, and
  Tyreus]{kwakernaak_1974_linear}
H.~Kwakernaak, R.~Sivan, and B.~N.~D. Tyreus.
\newblock \emph{Linear Optimal Control Systems}.
\newblock Wiley-Interscience, New York, 1974.

\bibitem[Lucchetti and Sans{\`o}(2020)]{lucchetti2020class}
R.~Lucchetti and F.~Sans{\`o}.
\newblock A class of sets where convergence in {H}ausdorff sense and in measure
  coincide.
\newblock \emph{Atti della Accademia Peloritana dei Pericolanti-Classe di
  Scienze Fisiche, Matematiche e Naturali}, 98\penalty0 (S2):\penalty0 9, 2020.

\bibitem[Munkres(2000)]{munkrestopology}
J.~Munkres.
\newblock \emph{Topology}.
\newblock Pearson College Div, 2000.

\bibitem[Platzman(1980)]{platzman_1980}
L.~K. Platzman.
\newblock Optimal infinite-horizon undiscounted control of finite probabilistic
  systems.
\newblock \emph{SIAM Journal on Control and Optimization}, 18\penalty0
  (4):\penalty0 362--380, 1980.

\bibitem[Rhenius(1974)]{rhenius_1974}
D.~Rhenius.
\newblock Incomplete information in {Markovian} decision models.
\newblock \emph{The Annals of Statistics}, 2\penalty0 (6):\penalty0 1327--1334,
  1974.

\bibitem[Rockafellar and Wets(2009)]{rockafellar_variational_2009}
R.~T. Rockafellar and R.~J.-B. Wets.
\newblock \emph{Variational Analysis}, volume 317.
\newblock Springer Science \& Business Media, Berlin, 2009.

\bibitem[Rudin(1987)]{walter_rudin_1966}
W.~Rudin.
\newblock \emph{{Real} and {Complex} {Analysis}}.
\newblock McGraw-Hill Book Company, Inc., New York, 1987.

\bibitem[Runggaldier and Stettner(1994)]{runggaldier_1994_approximations}
W.~J. Runggaldier and L.~Stettner.
\newblock \emph{Approximations of Discrete Time Partially Observed Control
  Problems}.
\newblock Applied Mathematics Monographs CNR, Giardini Editori, Pisa, 1994.

\bibitem[Sawaragi and Yoshikawa(1970)]{sawaragi_1970}
Y.~Sawaragi and T.~Yoshikawa.
\newblock Discrete-time {Markovian} decision processes with incomplete state
  observation.
\newblock \emph{The Annals of Mathematical Statistics}, 41\penalty0
  (1):\penalty0 78--86, 1970.

\bibitem[Sch{\"a}l(1972)]{schal_1972_continuous}
M.~Sch{\"a}l.
\newblock On continuous dynamic programming with discrete time parameter.
\newblock \emph{Z. Wahrscheinlichkeitstheorie verw. Gebiete}, 21:\penalty0
  279--288, 1972.

\bibitem[Sch{\"a}l(1993)]{schal_1993_average}
M.~Sch{\"a}l.
\newblock Average optimality in dynamic programming with general state space.
\newblock \emph{Mathematics of Operations Research}, 18\penalty0 (1):\penalty0
  163--172, 1993.

\bibitem[Shiryaev(1964 (in Russian))]{shiryaev_1966}
A.~Shiryaev.
\newblock On the theory of decision functions and control of a process of
  observation based on incomplete information.
\newblock pages 678--681, 1964 (in Russian).
\newblock Engl. transl. in \textit{Select. Transl. in Math., Statist. Probab.}
  6(1966), 162--188.

\bibitem[Smallwood and Sondik(1978)]{smallwood_1973_optimal}
R.~D. Smallwood and E.~J. Sondik.
\newblock The optimal control of partially observable {Markov} processes over a
  finite horizon.
\newblock \emph{Operations Research}, 21\penalty0 (5):\penalty0 1071--1088,
  1978.

\bibitem[Sondik(1978)]{sondik_1978_optimal}
E.~J. Sondik.
\newblock The optimal control of partially observable {Markov} processes over
  the infinite horizon: discounted costs.
\newblock \emph{Operations Research}, 26\penalty0 (2):\penalty0 282--304, 1978.

\bibitem[Wentzell(1996)]{wentzell_1996_course}
A.~D. Wentzell.
\newblock \emph{A Course in the Theory of Stochastic Processes}.
\newblock McGraw-Hill, New York, 1996.

\bibitem[Yushkevich(1976)]{yushkevich_1976_reduction}
A.~A. Yushkevich.
\newblock Reduction of a controlled {Markov} model with incomplete data to a
  problem with complete information in the case of {Borel} state and control
  spaces.
\newblock \emph{Theory of Probability and Its Applications}, 21\penalty0
  (1):\penalty0 153--158, 1976.

\end{thebibliography}

\end{document}